\documentclass[12pt,a4,paper]{article}
\usepackage{latexsym}
\usepackage{graphics}
\usepackage{amsthm}
\usepackage{amsmath}
\usepackage[symbol]{footmisc}
\DefineFNsymbols*{stars}{*{**}}
\setfnsymbol{stars}

\font\Bbb=msbm10 scaled\magstep1 \font\scriptBbb=msbm10
\font\scriptscriptBbb=msbm7\newfam
\Bbbfam\textfont\Bbbfam=\Bbb\scriptfont\Bbbfam=\scriptBbb
\scriptscriptfont\Bbbfam=\scriptscriptBbb

\def\CC{{\fam=\Bbbfam C}}

\def\KK{{\fam=\Bbbfam K}}

\def\NN{{\fam=\Bbbfam N}}
\def\RR{{\fam=\Bbbfam R}}
\def\Id{{\rm Id}}

\newtheorem{theorem}{Theorem}
\newtheorem{prop}[theorem]{Proposition}
\newtheorem{corollary}[theorem]{Corollary}
\newtheorem{lemma}[theorem]{Lemma}

\def\be{\begin{enumerate}}
\def\ee{\end{enumerate}}

\begin{document}

\parindent0pt

\title{On the maximal perimeter of sections of the cube}
\author{Hermann K\"onig (Kiel)\footnote{Part of the work was done when
the first-named author visited the University of Missouri-Columbia as a
Miller Distinguished Scholar}, Alexander Koldobsky \footnote{Partially
supported by the NSF grant DMS-1700036} (Columbia)}
\date{}

\maketitle

\begin{abstract}
We prove that the $(n-2)$-dimensional surface area (perimeter) of central hyperplane sections
of the $n$-dimensional unit cube is maximal for the hyperplane perpendicular to the vector
$(1,1,0,\dots,0)$. This gives a positive answer to a question of Pe{\l}czy\'nski who solved
the three dimensional case. We study both the real and the complex versions of this problem.
We also use our result to show that the answer to an analogue of the Busemann-Petty problem
for the surface area is negative in dimensions 14 and higher.
\end{abstract}


\section{Introduction, volume formulas and results}

A remarkable result of Ball \cite{B1} states that the hyperplane section of the $n$-cube $B_\infty^n$ perpendicular to $a_{max}:=\frac 1 {\sqrt 2}(1,1,0,\dots,0)$ has the maximal $(n-1)$-dimensional volume among all hyperplane sections, i.e. for any $a \in S^{n-1} \subset \RR^n$
$$vol_{n-1}(B_\infty^n \cap a^\perp) \le vol_{n-1}(B_\infty^n \cap a_{max}^\perp)=\sqrt 2,$$
where $a^\bot$ is the central hyperplane orthogonal to $a.$
Oleszkiewicz and Pe{\l}czy\'nski \cite{OP} proved the complex analogue of this result, with the same hyperplane $a_{max}^\perp$. \\

Pe{\l}czy\'nski \cite{P} asked whether the same hyperplane section is also maximal for intersections with the {\em boundary} of the $n$-cube, i.e. whether for all $a \in S^{n-1} \subset \RR^n$
$$ vol_{n-2}(\partial B_\infty^n \cap a^\perp) \le vol_{n-2}(\partial B_\infty^n \cap a_{max}^\perp)=2((n-2)\sqrt 2+1). $$
He proved it for $n=3$ when $vol_1(\partial B_\infty^3 \cap a^\perp)$ is the {\em perimeter} of the quadrangle or hexagon of intersection. In this paper, we answer Pe{\l}czy\'nski's question affirmatively for all $n \ge 3$. We also solve the complex version of this problem. For simplicity, we continue to call the quantity $vol_{n-2}(\partial B_\infty^n \cap a^\perp)$ the perimeter of the cubic section.\\

Ball \cite{B2} used his result to prove that the answer to the Busemann-Petty problem is negative in dimensions
10 and higher. The Busemann-Petty problem asks the following question. Suppose that origin-symmetric
convex bodies $K,L$ in $\RR^n$ satisfy
$$vol_{n-1}(K \cap a^\perp)\le vol_{n-1}(L \cap a^\perp)$$
for all $a\in S^{n-1}.$ Does it follow that the $n$-dimensional volume of $K$ is smaller than that of $L,$
i.e. $vol_n K\le vol_n L?$ The problem was solved as the result of work of many mathematicians, and the answer is affirmative for $n\le 4,$ and it is negative for $n\ge 5;$ see \cite{K} for details. Ball's result was one of
the steps of the solution. He showed that the answer is negative when $n\ge 10,$ $K$ is the unit cube and $L$ is the Euclidean ball of certain size in $\RR^n.$\\

We consider the following analogue of the Busemann-Petty problem for the surface area. Suppose that
origin-symmetric convex bodies $K,L$ in $\RR^n$ satisfy
$$vol_{n-2}(\partial K \cap a^\perp)\le vol_{n-2}(\partial L \cap a^\perp)$$
for all $a\in S^{n-1},$ i.e. the surface area (perimeter) of every central hyperplane section of $K$ is smaller than the same for $L.$ Does it follow that the surface area of $K$ is smaller than that of $L,$
i.e. $vol_{n-1}(\partial K)\le vol_{n-1}(\partial L)?$ We prove in Section \ref{BP} that the answer is negative for $n\ge 14$
and higher, with $K$ being the unit cube and $L$ the Euclidean ball of appropriate size in $\RR^n.$\\

To formulate our results precisely, let us introduce the following notations. Let $\KK \in \{\RR, \CC\}$, $\alpha= \frac 1 2$ for $\KK=\RR$ and $\alpha= \frac 1 {\sqrt \pi}$ for $\KK=\CC$. Let $|| \cdot ||_\infty$ and $| \cdot |$ denote the maximum and the Euclidean norm on $\KK^n$, respectively. Then
$$ B_\infty^n := \{ x \in \KK^n \ | \ || x ||_\infty \le \alpha \} $$
is the $n$-cube of volume 1 in $\KK^n$. For $\KK = \CC$, we identify $\KK^n = \RR^{2n}$ for volume calculations, i.e. we consider $vol_{2n}$ and $vol_{2n-2}$ for the polydisc $B_\infty^n$ and its complex hyperplane sections, respectively. For $a \in \KK^n$ with $|a|=1$ and $t \in \KK$, the {\em parallel section function} $A$ is defined by
$$ A_{l(n-1)}(a,t) := vol_{l(n-1)}(B_\infty^n \cap (a^\perp + \alpha t a)) , $$
where $l=1$ if $\KK=\RR$ and $l=2$ if $\KK=\CC$ and $a^\perp := \{ x \in \KK^n \ | \ \langle x , a \rangle =0 \}$. This gives the volume of the hyperplane section of the $n$-cube perpendicular to $a$ and at distance $\alpha t$ to the origin. We put $A_{l(n-1)}(a)=A_{l(n-1)}(a,0)$. Then Ball's result and Oleszkiewicz and
Pe{\l}czy\'nski's complex analogue state that for all $a \in \KK^n$ with $|a|=1$ we have
$$ A_{l(n-1)}(a) \le A_{l(n-1)}(a_{max}) = (\sqrt 2) ^l .$$
The lower bound $1 = A_{l(n-1)}(a_{min}) \le A_{l(n-1)}(a)$, $a_{min}=(1,0,\cdots,0)$, was shown earlier by Hensley \cite{H}. \\

For $a \in \KK^n$ with $|a|=1$ we define the {\em perimeter} of the cubic section by $a^\perp$ as
$$ P_{l(n-2)}(a) := vol_{l(n-2)}(\partial B_\infty^n \cap a^\perp) , $$
with $l$ as before. The main result of this paper answers Pe{\l}czy\'nski's problem affirmatively:
\begin{theorem}\label{th1}
Let $n \ge 3$ and $a_{max} := \frac 1 {\sqrt 2} (1,1,0, \cdots ,0) \in \KK^n$. Then for any $a \in \KK^n$ with $|a|=1$ we have
\begin{equation}\label{eq1}
P_{l(n-2)}(a) \le P_{l(n-2)}(a_{max}) ,
\end{equation}
where $l=1$ if $\KK=\RR$ and $l=2$ if $\KK=\CC$. We have \\
$$P_{n-2}(a_{max}) = 2 ( (n-2) \sqrt 2 +1 ) \quad , \quad \KK=\RR$$
and
$$P_{2(n-2)}(a_{max}) = 2 \pi ( (n-2) 2 +1 ) \quad , \quad \KK=\CC.$$
\end{theorem}

\vspace{0,4cm}
For $a \in \KK^n$ with $|a|=1$ let $a^\star$ denote the non-increasing rearrangement of the sequence $(|a_k|)_{k=1}^n$. Since the volume is invariant under coordinate permutations and sign changes, which in the complex case means rotations of coordinate discs, we have
$$A_{l(n-1)}(a,t) = A_{l(n-1)}(a^\star,|t|)$$
and
$$P_{l(n-2)}(a) = P_{l(n-2)}(a^\star).$$
Therefore, we will {\em generally} assume in this paper that $a = (a_k)_{k=1}^n$ satisfies $a_1 \ge \cdots \ge a_n \ge 0$ and $|a|=1$ as well as $t \ge 0$. For the parallel section function, the following formulas hold
\begin{equation}\label{eq2}
A_{n-1}(a,t) = \frac 2 \pi \; \int \limits_0^\infty \; \prod\limits_{k=1}^n \; \frac{\sin(a_ks)}{a_ks} \cos (ts) \ ds \quad , \quad \KK = \RR ,
\end{equation}
\begin{equation}\label{eq3}
A_{2(n-1)}(a,t)= \frac 1 2 \; \int\limits_0^\infty\; \prod\limits_{k=1}^n \; j_1 (a_ks) \; J_0(ts) \ s \ ds  \quad , \quad \KK=\CC ,
\end{equation}
where $j_1(t) = 2 \frac {J_1(t)} t$ and $J_\nu$ denote the Bessel functions of index $\nu$. If $a_k=0$, $\frac {\sin (a_ks)}{a_ks}$ and $j_1(a_ks)$ have to be read as 1 in formulas \eqref{eq2} and \eqref{eq3}. Formula \eqref{eq2} was shown already by P\'olya \cite{Po} in 1913, and was used by Ball \cite{B1} in the proof of his result. Both formulas can be proven by taking the Fourier transform of $A_{l(n-1)}(a, \cdot)$, using Fubini's theorem and taking the inverse Fourier transform, cf. e.g. Koldobsky, Theorem 3.1 \cite{K} or K\"onig, Koldobsky \cite{KK1} and \cite{KK2}. The $\frac {\sin t} t$ and $j_1(t)$ functions occur as Fourier transforms of the interval in $\RR$ and the disc in $\CC=\RR^2$, respectively. For the complex case cf. also Oleszkiewicz, Pe{\l}czy\'nski \cite{OP}. To prove Theorem \ref{th1}, we use the following formulas for the perimeter.

\begin{prop}\label{prop2}
For any $a = (a_k)_{k=1}^n \in S^{n-1} \subset \RR^n$
\begin{equation}\label{eq4}
P_{n-2}(a)= 2 \  \sum_{k=1}^n \sqrt {1-a_k^2} \; \frac 2 \pi \; \int_0^\infty  \prod_{j=1,j \ne k}^n \frac {\sin(a_js)}{a_js} \ \cos(a_ks) \ ds \quad , \quad \KK=\RR ,
\end{equation}
\begin{equation}\label{eq5}
P_{2(n-2)}(a)= 2 \pi \  \sum_{k=1}^n (1-a_k^2) \; \frac 1 2 \; \int_0^\infty  \prod_{j=1,j \ne k}^n j_1(a_js) \ J_0(a_ks) \ s \ ds \quad , \quad \KK=\CC .
\end{equation}
\end{prop}

In Ball's result, the integral in \eqref{eq2} for $t=0$ is estimated by using H\"older's inequality if $a_1 \le \frac 1 {\sqrt 2}$, which is natural since in the extremal case ($a_1=a_2=\frac 1 {\sqrt 2}, a_j=0, j>3$) the integrand is non-negative. In \eqref{eq4} and \eqref{eq5} we have weighted sums of integrals where the integrands are non-positive in the extremal case. Therefore, estimating $P_{l(n-2)(a)}$ requires further methods in addition to Ball's techniques and inequalities or those of Oleszkiewicz and Pe{\l}czy\'nski. One idea is to consider the perimeter estimate as a constrained optimization problem, in view of the following two results. We continue to denote $l=1$ if $\KK=\RR$ and $l=2$ if $\KK=\CC$.

\begin{prop}\label{prop3}
For any $a = (a_k)_{k=1}^n \in S^{n-1} \subset \RR^n$ and $k \in \{1, \cdots , n\}$, define
\begin{eqnarray}\label{eq6}
D_k(a) :=  \begin{cases}
                  & \frac 2 \pi \; \int_0^\infty  \prod_{j=1,j \ne k}^n \frac {\sin(a_js)}{a_js} \ \cos(a_ks) \ ds \quad , \quad \KK=\RR \\
                  & \frac 1 2 \; \int_0^\infty  \prod_{j=1,j \ne k}^n j_1(a_js) \ J_0(a_ks) \ s \ ds \; \ , \quad \KK=\CC
                  \end{cases}
                  \ \Bigg\} \ .
\end{eqnarray}
Then
\begin{equation}\label{eq7}
\sum_{k=1}^n \ D_k(a) = (n-1) \ A_{l(n-1)}(a) .
\end{equation}
\end{prop}

\begin{prop}\label{prop4}
For any $a = (a_k)_{k=1}^n \in S^{n-1} \subset \RR^n$ and $k \in \{1, \cdots , n\}$,
\begin{equation}\label{eq8}
D_k(a) \le A_{l(n-1)}(a) .
\end{equation}
\end{prop}

The proof of Proposition \ref{prop4} also yields the following estimate for the parallel section function

\begin{corollary}\label{cor5}
For any $a \in \KK^n$ with $|a|=1$ and $t>0$ we have
$$A_{l(n-1)}(a,t) \le \sqrt { \frac 2 {1+t^2} }^l .$$
\end{corollary}

Ball's proof relies on the non-trivial estimate $f(p) \le f(2) = 1$ for the function
$$f(p) := \sqrt{ \frac p 2 } \ \frac 2 \pi \ \int_0^\infty  \left|\frac {\sin(t)} t \right|^p \ dt \ , $$
since then in the real case for all $0 < a_n \le \cdots \le a_1 \le \frac 1 {\sqrt 2}$ with $\sum_{k=1}^n a_k^2 =1$ we find by using H\"older's inequality with $p_k:=a_k^{-2} \ge 2$
\begin{align}\label{eq9}
A_{n-1}(a) & \le \prod_{k=1}^n (\frac 2 \pi \int_0^\infty \left|\frac {\sin(a_ks)}{a_ks}\right|^{a_k^{-2}} \ ds )^{a_k^2} \nonumber \\
& = (\prod_{k=1}^n f(a_k^{-2}))^{a_k^2} \sqrt 2 \le \sqrt 2 .
\end{align}

The constrained approximation approach suffices to prove Theorem \ref{th1}, except when, in the real case, $a_1$ is in a small interval around $\frac 1 {\sqrt 2}$. To prove Theorem \ref{th1} also in this case, we need additional information on the function $f$:

\begin{prop}\label{prop6}
Define $f:(1,\infty) \to \RR_+$ by
$$f(p):= \sqrt { \frac p 2 } \ \frac 2 \pi \ \int_0^\infty \left|\frac {\sin(t)}  t \right|^p \ dt. $$
Then \\
(a) $\lim_{p \to \infty} f(p) = \sqrt{\frac 3 \pi}$  and  $f(\frac 9 4) < \sqrt{\frac 3 \pi}$ . \\
(b) $f(\sqrt 2 + \frac 1 2) < \frac {51}{50}$ . \\
(c) $f|_{[\sqrt 2 + \frac 1 2 , \frac 9 4]}$ is decreasing and convex.
\end{prop}

\begin{prop}\label{prop7}
For all $p \ge \frac 9 4$, $f(p) \le \sqrt{\frac 3 \pi }$.
\end{prop}

Using the convexity of $f$ and the estimates for $f(p)$ for $p=\sqrt 2 +1/2$ and $p= \frac 9 4$, we may improve the general estimate \eqref{eq9} for certain sequences $a$ with $a_1$ close to $\frac 1 {\sqrt 2}$, which will essentially suffice to prove Theorem \ref{th1} in these cases. This works since $f$ is strictly smaller near $\infty$ than in $p=2$ where $f(2)=1$. In the complex case, the function $f$ is replaced by
$$\tilde{f}(p):= \frac p 2 \frac 1 2 \int_0^\infty |j_1(s)|^p s \ ds \ , $$
where also $\tilde{f}(p) \le \tilde{f}(2)=1$ for all $p \ge 2$, cf.  Oleszkiewicz, Pe{\l}czy\'nski \cite{OP}. However, in this case $\lim_{p \to \infty} \tilde{f}(p) = \tilde{f}(2) = 1$. Therefore, no analogue of Proposition \ref{prop6} (a), (b) and Proposition \ref{prop7} is possible in the complex case. Fortunately, in the complex case, the perimeter formula given by \eqref{eq5} is easier to estimate since it does not contain a square root in the weights of the integrals, and the constraint technique works for all sequences $a$. \\

\section{Constrained optimization}

We start by proving the formulas for the perimeter. \\

{\bf Proof of Proposition \ref{prop2}}. \\
Let $a=(a_k)_{k=1}^n \in S^{n-1}$, $a_1 \ge \cdots a_n \ge 0$ and $x \in \KK^n$. We write $a=(a_1,\tilde{a})$, $x=(x_1,\tilde{x})$ with $\tilde{a}=(a_k)_{k=2}^n, \tilde{x}=(x_k)_{k=2}^n \in \KK^{n-1}$, a notation also used in the following proofs. In the real case $\KK=\RR$, the $(n-1)$-dimensional hyperplane $a^\perp$ intersects the boundary $\partial B_\infty^n$ in $2n$ $(n-2)$-dimensional (typically non-central) sections of an $(n-1)$-cube, namely for $x_j = \pm \frac 1 2$, $j=1, \cdots , n$. For $x_1 = - \frac 1 2$ we need to calculate
$$ vol_{n-2} \{ \tilde{x} \in \RR^{n-1} \ | \ \langle \tilde{x} , \tilde{a} \rangle = \frac 1 2 a_1 \} . $$
Let $a_j' := \frac {a_j}{\sqrt{1-a_1^2}}$, $j=1, \cdots , n$ and $\tilde{a}' := (a_j')_{j=2}^n$. Then $|\tilde{a}'|^2 = \sum_{j=2}^n a_j'^2 = 1$. Using \eqref{eq2}, we find
\begin{align*}
vol_{n-2} & \{ \ \tilde{x} \in \RR^{n-1} \ | \ \langle \tilde{x} , \tilde{a} \rangle = \frac 1 2 a_1 \ \} = A_{n-2}(\tilde{a}',a_1') \\
& = \frac 2 \pi \int_0^\infty \prod_{j=2}^n \frac {\sin(a_j'r)}{a_j'r} \cos(a_1'r) dr =
\sqrt{1-a_1^2} \ \frac 2 \pi \int_0^\infty \prod_{j=2}^n \frac {\sin(a_js)}{a_js} \ \cos(a_1s) ds \ .
\end{align*}
The same holds for $x_1 = + \frac 1 2$ and similarly for $x_j= \pm \frac 1 2$, so that
$$P_{n-2}(a)= 2 \  \sum_{k=1}^n \sqrt {1-a_k^2} \; \frac 2 \pi \; \int_0^\infty  \prod_{j=1,j \ne k}^n \frac {\sin(a_js)}{a_js} \ \cos(a_ks) \ ds \ , $$
which proves \eqref{eq4}.\\
In the complex case $\KK=\CC$, we have to consider the intersection of $a^\perp$ with $x_j=\frac 1 {\sqrt \pi} \exp(i \theta)$ for all $\theta \in [0,\pi)$, and use \eqref{eq3} instead of \eqref{eq2}. Then
\begin{align*}
P_{2(n-2)}(a) &= 2 \pi \  \sum_{k=1}^n \frac 1 2 \; \int_0^\infty  \prod_{j=1,j \ne k}^n j_1(a_j'r) \ J_0(a_k'r) \ r \ dr \\
&= 2 \pi \  \sum_{k=1}^n (1-a_k^2) \; \frac 1 2 \; \int_0^\infty  \prod_{j=1,j \ne k}^n j_1(a_js) \ J_0(a_ks) \ s \ ds \ ,
\end{align*}
which yields formula \eqref{eq5}. \hfill $\Box$ \\

Pe{\l}czy\'nski \cite{P} proved Theorem \ref{th1} for $n=3$ in the real case by considering three affine independent points on the boundary of the cube and their antipodals, calculating the perimeter of the (possibly non-planar) hexagon defined that way. This perimeter then turned out to be maximal in the case that the hexagon degenerates into a rectangle perpendicular to e.g. $(1,1,0)$, which is planar. We give the easy direct proof of Theorem \ref{th1} for $n=3$, $\KK=\RR$ by using Proposition \ref{prop2}. \\

{\bf Proof of Theorem 1} for $n=3, \KK=\RR$: \\
Let $a_1 \ge a_2 \ge a_3 \ge 0$, $a_1^2+a_2^2+a_3^2 = 1$. Calculating the integrals in \eqref{eq4}, we find that
\begin{eqnarray*}
\frac 1 2 P_1(a) = \begin{cases}
                   \quad \quad \quad \frac 1 {a_1} ( \sqrt{1-a_2^2} + \sqrt{1-a_3^2} )   & , \quad a_1 \ge a_2 + a_3 \\
                   \sqrt{1-a_1^2}\frac{a_2+a_3-a_1}{2a_2a_3}+\sqrt{1-a_2^2}\frac{a_1+a_3-a_2}{2a_1a_3}+\sqrt{1-a_3^2}\frac{a_1+a_2-a_3}{2a_1a_2} & , \quad a_1 < a_2+a_3
                  \end{cases}
                  \Bigg\} \ .
\end{eqnarray*}
In the first case, the hyperplane intersects the cube in a rectangle, in the second case in a hexagon. \\

i) Assume first that $a_1 \ge a_2+a_3$. Then $(a_1-a_2)^2 \ge a_3^2 = 1 - a_1^2-a_2^2$, $1-a_2^2 \le a_1^2 + (a_1-a_2)^2$ and
$$\frac 1 2 P_1(a) \le \frac 1 {a_1}( \sqrt{a_1^2+(a_1-a_2)^2} + \sqrt{a_1^2+a_2^2} ) = \sqrt{1-(1-x)^2} + \sqrt{1-x^2},$$
where $0 \le x:= \frac {a_2}{a_1} \le 1$. The right side is maximal for $x=0$ or $x=1$ with $\frac 1 2 P_1(a) \le \sqrt 2 + 1$. \\

ii) If $a_1 < a_2+a_3$, assume first that $a_1=a_2 \ge a_3 \ge0$. Then $\frac 1 {\sqrt 3} \le a_2 = a_1 \le \frac 1 {\sqrt 2}$, $1-a_3^2 = 2 a_1^2$ and, as easily seen by the above formula,
$$\frac 1 2 P_1(a) = \sqrt 2 + \frac 1 {a_1} ( \sqrt{1-a_1^2} - \sqrt{\frac 1 2 - a_1^2} ) \le \sqrt 2 +1 .$$

If $a_1 < a_2+a_3$, but $a_1 > a_2$, define $x$ by $a_3 = x (a_1-a_2)$ so that $x \ge 1$. Then
$$\frac 1 2 P_1(a) = \frac 1 2 (\sqrt{1-a_1^2}\ \frac{x-1} {xa_2}+\sqrt{1-a_2^2}\ \frac{x+1} {xa_1}+\sqrt{1-x^2(a_1-a_2)^2}\ \frac{a_1+a_2-x(a_1-a_2)}{a_1a_2}) =: \frac 1 2 \psi(x). $$
We have
$$\psi'(x)=(\frac{\sqrt{1-a_1^2}}{a_2}-\frac{\sqrt{1-a_2^2}}{a_1}) \frac 1 {x^2}-\sqrt{1-a_3^2}\ \frac{a_1-a_2}{a_1a_2} -\frac{x(a_1-a_2)^2}{\sqrt{1-a_3^2}} \frac{a_1+a_2-x(a_1-a_2)}{a_1a_2} . $$
If the factor of $\frac 1 {x^2}$ is negative, all summands are negative and $\psi'(x) \le 0$. If the factor is positive,
\begin{align*}
\psi'(x) &\le (\frac{\sqrt{1-a_1^2}}{a_2}-\frac{\sqrt{1-a_2^2}}{a_1}) - \sqrt{a_1^2+a_2^2} \ \frac{a_1-a_2}{a_1a_2} \\
&= \frac 1 {a_1a_2} (a_1 \sqrt{1-a_1^2}-a_2 \sqrt{1-a_2^2}-(a_1-a_2) \sqrt{a_1^2+a_2^2}) .
\end{align*}
This is negative as well: $\phi(y) := y \sqrt{1-y^2}$ satisfies $\phi'(y) = \frac{1-2y^2}{\sqrt{1-y^2}}$, so that
$a_1 \sqrt{1-a_1^2}-a_2 \sqrt{1-a_2^2} = (a_1-a_2) \frac{1-2y^2}{\sqrt{1-y^2}}$ for some $a_2<y<a_1$. But
$$\frac{1-2y^2}{\sqrt{1-y^2}} \le \frac{1-2a_2^2}{\sqrt{1-a_2^2}} \le \sqrt{1-a_2^2} = \sqrt{a_1^2+a_3^2} \le \sqrt{a_1^2+a_2^2} \ . $$
Hence $\psi'(x) \le 0$, so that $\psi(x) \le \psi(1)$ since $x \ge 1$. Therefore
$$\frac 1 2 P_1(a) \le \frac 1 {a_1} ( \sqrt{1-a_2^2} + \sqrt{1-a_3^2} \ ) , $$
and, since $x=1$, $a_1^2+a_2^2+(a_1-a_2)^2=1$, $a_2 = \frac 1 2 (a_1+\sqrt{2-3a_1^2})$, yielding
$$1-a_2^2 = \frac 1 2 (1+a_1^2) - \frac 1 2 a_1 \sqrt{2-3a_1^2} =: \phi_-(a_1) \ , \ 1-a_3^2 = a_1^2+a_2^2 = \frac 1 2 (1+a_1^2) + \frac 1 2 a_1 \sqrt{2-3a_1^2} =: \phi_+(a_1) , $$
with $\frac 1 {\sqrt 2} \le a_1 \le \sqrt{\frac 2 3}$ so that
$$\frac 1 2 P_1(a) \le \frac 1 {a_1} (\sqrt{\phi_-(a_1)} + \sqrt{\phi_+(a_1)}) \le \sqrt 2 + 1 , $$
the maximal value being attained for $a_1 = \frac 1 {\sqrt 2}$. \hfill $\Box$ \\

For $\KK=\RR$ and $n=4$, integration of formula \eqref{eq4} yields three cases

\begin{eqnarray}\label{eq10}
\frac 1 2 P_2(a) := \begin{cases}
                  & \sqrt{1-a_1^2} \ \frac{a_2+a_3-a_1}{2a_2a_3} + \sqrt{1-a_2^2} \ \frac{a_1+a_3-a_2}{2a_1a_3} + \sqrt{1-a_3^2} \ \frac{a_1+a_2-a_3}{2a_1a_2} \\
                  & \quad \quad \quad + \sqrt{1-a_4^2} \ (\frac 1 {a_1} - \frac{a_4^2(a_2+a_3-a_1)^2}{4a_1a_2a_3}) \quad \ , \ a_1 < a_2+a_3-a_4 \\ \\

                  & - \frac{(a_2+a_3+a_4-a_1)^2}{8a_1a_2a_3a_4} \ (-a_1 \sqrt{1-a_1^2} + a_2  \sqrt{1-a_2^2} + a_3 \sqrt{1-a_3^2} + a_4 \sqrt{1-a_4^2}) \\
                  & \quad + \frac 1 {a_1} (\sqrt{1-a_2^2} + \sqrt{1-a_3^2} + \sqrt{1-a_4^2}) \ , \ a_2+a_3-a_4 \le a_1 \le a_2+a_3+a_4 \\ \\

                  & \frac 1 {a_1} (\sqrt{1-a_2^2} + \sqrt{1-a_3^2} + \sqrt{1-a_4^2}) \quad \quad \ , \ a_2+a_3+a_4 < a_1
                  \end{cases} \ .
\end{eqnarray}
These formulas can be derived e.g. by using formula (2.1) of K\"onig, Koldobsky \cite{KK1}. \\

{\bf Proof of Proposition \ref{prop3}.}\\
We first give a geometric proof in the real case. \\
a) Let $a_1 \ge \cdots \ge a_n \ge 0$, $\sum_{k=1}^n a_k^2 =1$, $a_k':= \frac {a_k}{\sqrt{1-a_1^2}}$, $\tilde{a}' := (a_k')_{k=2}^n$ and $D_j(a)$ be given as in \eqref{eq6}. By transformation of variables
$$D_1(a) = \frac 1 {\sqrt{1-a_1^2}} \ \frac 2 \pi \int_0^\infty \prod_{j=2}^n \frac{\sin(a_j'r)}{a_j'r} \ \cos(a_1'r) \ dr = \frac 1 {\sqrt{1-a_1^2}} \ A_{n-2}(\tilde{a}',a_1') ,$$
in terms of the $(n-2)$-dimensional volume of the section of $B_\infty^{n-1}$ perpendicular to $\tilde{a}'$ and at distance $\frac 1 2 a_1'$ to the origin of $B_\infty^{n-1}$. Since $\frac 1 2 \frac 1 {\sqrt{1-a_1^2}}$ is the height of the $(n-1)$-dimensional pyramid $P(1)$ with vertex in 0 and base being the above $(n-2)$-dimensional section,
$$vol_{n-1}(P(1))= \frac 1 {n-1} \ A_{n-2}(\tilde{a}',a_1') \ \frac 1 2 \frac 1 {\sqrt{1-a_1^2}} \ . $$
Summing up the volumes of these pyramids, also for opposite sections, yields
$$A_{n-1}(a) = 2 \sum_{k=1}^n vol_{n-1}(P(k)) = \frac 1 {n-1} \sum_{k=1}^n D_k(a) , $$
which is \eqref{eq7}. \\

b) We now give a second, analytic proof of \eqref{eq7}, based on integration by parts, using
$$ \frac d {ds} (\frac {\sin(a_js)}{a_js}) = \frac 1 s ( \cos(a_js) - \frac{\sin(a_js)}{a_js} ) , $$
if all $a_j$ are $>0$. Then
\begin{align*}
D_1(a) &= \left[\frac 2 \pi \prod_{j=2}^n \frac{\sin(a_js)}{a_js} \frac{\sin(a_1s)}{a_1}\right]_{s=0}^\infty -
\frac 2 \pi \int_0^\infty \sum_{j=2}^n \prod_{k=2,k \ne j}^n \frac{\sin(a_ks)}{a_ks} \left(\frac{\cos(a_js)}{s} - \frac{\sin(a_js)}{a_js^2}\right) \frac{\sin(a_1s)}{a_1} ds \\
& = (n-1) A_{n-1}(a) - \sum_{j=2}^n D_j(a) \ ,
\end{align*}
so that $\sum_{j=1}^n D_j(a) = (n-1) A_{n-1}(a)$. If some $a_k$ are zero, the corresponding $D_k(a)$ equals $A_{n-1}(a)$, and \eqref{eq7} follows by integration by parts only for those $k$ where $a_k \ne 0$. \\

c) The integration by parts technique also works in the complex case $\KK=\CC$, using
$$\frac d {ds} j_1(s) = 2 \frac d {ds} \left(\frac {J_1(s)} s\right) = - 2 \frac{J_2(s)} s = 2 \frac{J_0(s)} s - 4 \frac{J_1(s)} {s^2} $$
and $\frac d {ds} (s J_1(s)) = s J_0(s)$. For these formulas on Bessel functions, cf. Watson \cite{W}.  \hfill $\Box$ \\

{\bf Proof of Proposition \ref{prop4}.}\\
We first consider the real case. To show $D_k(a) \le A_{n-1}(a)$, we may assume without loss of generality that $k=1$ and $a_1 >0$. We will not use any inequality between the coordinates of $a$ in this proof, but assume that $a_k \ge 0$ for all $k$. Again, let $a_j' = \frac{a_j}{\sqrt{1-a_1^2}}$ for $j=1, \cdots, n$,
$\tilde{a}' := (a_j')_{j=2}^n \in \RR^{n-1}$. Then $\sum_{j=2}^n a_j'^2 =1$ so that by transformation of variables and \eqref{eq6} in dimension $n-1$
\begin{align*}
D_1(a) &= \frac 1 {\sqrt{1-a_1^2}} \ \frac 2 \pi \int_0^\infty \prod_{j=2}^n \frac{\sin(a_j'r)}{a_j'r} \ \cos(a_1'r) \ dr \\
&= \frac 1 {\sqrt{1-a_1^2}} \ vol_{n-2} \{\ \tilde{x} \in B_\infty^{n-1} \ | \ \langle \tilde{x} , \tilde{a}' \rangle = \frac 1 2 a_1' \ \} .
\end{align*}
By Brunn-Minkowski, we have for any $t \in \RR$ with $|t| \le a_1'$
$$vol_{n-2} \{\tilde{x} \in B_\infty^{n-1} \ | \ \langle \tilde{x} , \tilde{a}' \rangle = \frac 1 2 a_1' \} \le
vol_{n-2} \{\tilde{x} \in B_\infty^{n-1} \ | \ \langle \tilde{x} , \tilde{a}' \rangle = \frac 1 2 t \} . $$
Therefore
$$D_1(a) \le \frac 1 {a_1'} \frac 1 {\sqrt{1-a_1^2}} \ vol_{n-1} \{\tilde{x} \in B_\infty^{n-1} \ | \ | \langle \tilde{x} , \tilde{a}' \rangle | \le \frac 1 2 a_1' \} . $$
Define $T : \RR^{n-1} \to \RR^n$ by $T(x) := (-\frac{\langle \tilde{x},\tilde{a}' \rangle}{a_1'} , \tilde{x})$. Then $T$ maps the slab \newline
$\{\tilde{x} \in B_\infty^{n-1} \ | \ | \langle \tilde{x} , \tilde{a}' \rangle | \le \frac 1 2 a_1' \}$ in dimension $n-1$ into a central section of $B_\infty^n$,
\begin{align*}
T \{\tilde{x} \in B_\infty^{n-1} \ | \ | \langle \tilde{x} , \tilde{a}' \rangle | \le \frac 1 2 a_1' \} & =
\{y=(x_1,\tilde{x}) \in B_\infty^n \ | \ x_1 a_1 + \langle \tilde{x} , \tilde{a}' \rangle = 0 \} \\
& = \{y \in B_\infty^n \ | \ \langle y, a \rangle = 0 \} .
\end{align*}
Recall here that we normalized $B_\infty^{n-1}$ and $B_\infty^n$ to have volume 1. Since
$$ T^*T = \Id + \frac 1 {a_1'^2} \ \langle \cdot , \tilde{a}'\rangle \ \tilde{a}' , $$
$T^*T$ has an $(n-2)$-fold eigenvalue 1 and one eigenvalue  $1 + \frac 1 {a_1'^2}$ (with eigenvector $\tilde{a}'$ of norm 1) so that
$$ \sqrt{\det(T^*T)} = \frac 1 {a_1'} \sqrt{a_1'^2 + 1} = \frac 1 {a_1'} \frac 1 {\sqrt{1-a_1^2}} .$$
Therefore
\begin{align*}
D_1(a) &\le \sqrt{ \det(T^*T) } \ vol_{n-1} \{ \tilde{x} \in B_\infty^{n-1} \ | \ | \langle \tilde{x} , \tilde{a}' \rangle | \le \frac 1 2 a_1' \} \\
& = vol_{n-1} \{y \in B_\infty^n \ | \ \langle y, a \rangle = 0 \} = A_{n-1}(a) .
\end{align*}

The complex case requires only minor modifications. In that case
\begin{align*}
D_1(a) &= \frac 1 {1-a_1^2} \ vol_{2(n-1)} \{ \tilde{x} \in B_\infty^{n-1} \ | \ \langle \tilde{x} , \tilde{a}' \rangle = \frac 1 {\sqrt{\pi}} a_1' \} \\
& \le \frac 1 {a_1'^2} \frac 1 {1-a_1^2} \ vol_{2(n-1)} \{ \tilde{x} \in B_\infty^{n-1} \ | \ | \langle \tilde{x} , \tilde{a}' \rangle | \le \frac 1 {\sqrt{\pi}} a_1' \} .
\end{align*}
Define $T: \CC^{n-1} \to \CC^n$ also by $T(x) := (-\frac{\langle \tilde{x},\tilde{a}' \rangle}{a_1'} , \tilde{x})$, mapping the slab in $\CC^{n-1}$ into the central section in $\CC^n$ defined by $a^\perp$. In the complex case $\sqrt{\det(T^*T)} = \frac 1 {a_1'^2} (1+a_1'^2) = \frac 1 {a_1'^2} \frac 1 {1-a_1^2} $, so that
$$D_1(a) \le vol_{2(n-1)} \{ y \in B_\infty^n \ | \ \langle y , a \rangle = 0 \} = A_{2(n-1)}(a). $$  \hfill $\Box$ \\

{\bf Proof of Corollary \ref{cor5}. } \\
Let $a \in \KK^n$, $|a|=1$ and $t > 0$. Put $a_j' = \frac{a_j}{\sqrt{1+t^2}}$, $t':= \frac t {\sqrt{1+t^2}}$. Then $(a',t') \in \KK^{n+1}$, $|(a',t')|=1$. We get from \eqref{eq2} and \eqref{eq3} by transformation of variables
$A_{l(n-1)}(a,t) = \frac 1 {\sqrt{1+t^2}} \ D((a',t'))$, where e.g. in the real case
$$D((a',t')) = \frac 2 {\pi} \ \int_0^\infty \prod_{j=1}^n \frac{\sin(a_j's)}{a_j's} \ \cos(t's) \ ds . $$
By Proposition \ref{prop4}, applied to $(a',t') \in S^n \subset \RR^{n+1}$, $D((a',t')) \le A_n((a',t')) \le \sqrt 2$.
Similarly, for $\KK=\CC$, $D((a',t')) \le A_{2n}((a',t')) \le 2$. Therefore
$$A_{l(n-1)} \le \sqrt{\frac 2 {1+t^2}}^l . $$ \hfill $\Box$ \\

In the case that the largest coordinate $a_1$ of $a \in S^{n-1}$ satisfies $a_1 > \frac 1 {\sqrt 2}$, Ball \cite{B1} showed by projecting $a^\perp$ onto $a_{min}^\perp = (1,0,\cdots,0)^\perp$ that
\begin{equation}\label{eq11}
A_{n-1}(a) \le \frac 1 {a_1} \ A_{n-1}(a_{min}) = \frac 1 {a_1} .
\end{equation}
The complex analogue of this is, again if $a_1 > \frac 1 {\sqrt 2}$,
\begin{equation}\label{eq12}
A_{2(n-1)}(a) \le \frac 1 {a_1^2} \ A_{2(n-1)}(a_{min}) = \frac 1 {a_1^2} ,
\end{equation}
cf. Oleszkiewicz, Pe{\l}czy\'nski \cite{OP}. \\

We now prove Theorem \ref{th1}, except in the real case when $a_1 \in (\sqrt{\sqrt 2 -1},\frac 1 {\sqrt{\sqrt 2 + \frac 1 2}}) \simeq (0.643,0.723)$, i.e. when $a_1$ is close to $\frac 1 {\sqrt 2}$. This is done using the constraints given by Propositions \ref{prop3} and \ref{prop4}. \\

{\bf Proof of Theorem \ref{th1}.} \\
(a) We first verify the result in the complex case $\KK = \CC$ which is easier to prove. We have for $a_{max}$
$$D_1(a_{max})=D_2(a_{max})= \frac 1 2 \int_0^\infty j_1(\frac s {\sqrt2}) \ J_0(\frac s {\sqrt2}) \ s \ ds = 2 \int_0^\infty J_1(t) \ J_0(t) \ dt = 1 ,$$
$$D_j(a_{max})= \frac 1 2 \int_0^\infty j_1(\frac s {\sqrt2})^2 \ s \ ds = 4 \int_0^\infty J_1(t)^2 \frac {dt} t = 2 \ , \ j>2 .$$
For these integrals, cf. Gradstein, Ryshik \cite{GR} or Watson \cite{W}. Hence by \eqref{eq5}
\begin{equation}\label{eq13}
P_{2(n-1)}(a_{max}) = 2 \pi (1+(n-2) 2) = 2 \pi (2n-3) .
\end{equation}

Now consider $a=(a_k)_{k=1}^n \in S^{n-1}$ with $a_1 \ge \cdots \ge a_n \ge 0$. By Proposition \ref{prop2}
$$\frac 1 {2 \pi} P_{2(n-2)}(a) = \sum_{k=1}^n (1-a_k^2) \ D_k(a) ,$$
and using Propositions \ref{prop3} and \ref{prop4}, we have
$$\frac 1 {2 \pi} P_{2(n-2)}(a) \le \sup \{ \sum_{k=1}^n (1-a_k^2) \ C_k \ | \ 0 \le C_k \le A_{2(n-1)}(a) , \ \sum_{k=1}^n C_k = (n-1) \ A_{2(n-1)}(a) \} . $$
Since $(1-a_k^2)_{k=1}^n$ is increasing in $k$, the sum $\sum_{k=1}^n (1-a_k^2) \ C_k$ will be maximal under the given restrictions, if the sequence $(C_k)_{k=1}^n$ is increasing as well which, in fact, means that $C_1=0$, $C_2 =  \cdots =C_n = A_{2(n-1)}(a)$. Therefore
\begin{align*}
\frac 1 {2 \pi} P_{2(n-2)}(a) & \le \sum_{k=2}^n (1-a_k^2) \ A_{2(n-1)}(a) \\
& = [(n-1)- \sum_{k=2}^n a_k^2] \ A_{2(n-1)}(a) = (n-2+a_1^2) \ A_{2(n-1)}(a) .
\end{align*}
If $a_1 \le \frac 1 {\sqrt 2}$, we use that $A_{2(n-1)}(a) \le A_{2(n-1)}(a_{max}) = 2$ by \cite{OP}, so that with \eqref{eq13}
$$\frac 1 {2 \pi} P_{2(n-2)}(a) \le (n - \frac 3 2) \ 2 = \frac 1 {2 \pi} P_{2(n-2)}(a_{max}) . $$
If $a_1 > \frac 1 {\sqrt 2}$, we use that by \eqref{eq12} $A_{2(n-1)}(a) \le \frac 1 {a_1^2}$, so that
$$\frac 1 {2 \pi} P_{2(n-2)}(a) \le (n-2+a_1^2) \frac 1 {a_1^2} = \frac {n-2}{a_1^2} +1 \le (n-2) \ 2 + 1 = \frac 1 {2 \pi} P_{2(n-2)}(a_{max}) . $$
This proves Theorem \ref{th1} in the case of complex scalars. \\

(b) In the real case $\KK=\RR$, we have for $a_{max}$
$$D_1(a_{max})=D_2(a_{max}) = \frac 2 \pi \int_0^\infty \frac {\sin(\frac s {\sqrt 2})}{\frac s {\sqrt 2}} \ \cos(\frac s {\sqrt 2}) \ ds =
\frac 1 {\sqrt 2} \frac 2 \pi \int_0^\infty \frac{\sin(t)} t \ dt = \frac 1 {\sqrt 2} ,$$
$$D_j(a_{max}) = \frac 2 \pi \int_0^\infty \left(\frac{\sin(\frac s {\sqrt 2})}{\frac s {\sqrt 2}}\right)^2 \ ds  = \sqrt 2 \frac 2 \pi \int_0^\infty \left(\frac{\sin(t)} t\right)^2 \ dt = \sqrt 2\ , \ j>2 .$$
Hence by \eqref{eq4}
\begin{equation}\label{eq14}
P_{n-2}(a_{max}) = 2 ( (n-2) \sqrt 2 +1 ) .
\end{equation}
Now let $a=(a_k)_{k=1}^n \in S^{n-1}$ be arbitrary with $a_1 \ge \cdots \ge a_n \ge 0$. By Propositions \ref{prop2}, \ref{prop3} and \ref{prop4} we get, similarly as in part (a),
$$\frac 1 2 P_{n-2}(a) \le \sup \{ \sum_{k=1}^n \sqrt{1-a_k^2} \ C_k \ | \ 0 \le C_k \le A_{n-1}(a) , \ \sum_{k=1}^n C_k = (n-1) \ A_{n-1}(a) \} . $$
Since also $(\sqrt{1-a_k^2} \ )_{k=1}^n$ is increasing in $k$, the supremum is attained for increasing $C_k$ as well and, in fact, for $C_1=0$,
$C_2 = \cdots = C_k = A_{n-1}(a)$ so that
\begin{equation}\label{eq15}
\frac 1 2 P_{n-2}(a) \le \sum_{k=2}^n \sqrt{1-a_k^2} \ A_{n-1}(a) .
\end{equation}
Since $\phi(x) = \sqrt{1-x}$ is concave,
$$\frac 1 {n-1} \sum_{k=2}^n \phi(a_k^2) \le \phi(\frac 1 {n-1} \sum_{k=2}^n a_k^2) = \phi(\frac 1 {n-1} (1-a_1^2)) .$$
Hence
\begin{equation}\label{eq16}
\frac 1 2 P_{n-2}(a) \le (n-1) \sqrt{1- \frac{1-a_1^2}{n-1}} \ A_{n-1}(a)  \le (n-1- \frac{1-a_1^2} 2) \ A_{n-1}(a).
\end{equation}
If $a_1 \le \frac 1 {\sqrt 2}$, we use that $A_{n-1}(a) \le \sqrt 2$ by \cite{B1} to get
$$\frac 1 2 P_{n-2}(a) \le (n-2) \sqrt 2 + \frac 3 4 \sqrt 2 . $$
If $a_1 > \frac 1 {\sqrt 2}$, we use that $A_{n-1}(a) \le \frac 1 {a_1}$ by \eqref{eq11} and also find
$$\frac 1 2 P_{n-2}(a) \le (n- \frac 3 2 + a_1^2) \frac 1 {a_1} \le (n-2) \sqrt 2 + \frac 3 4 \sqrt 2 . $$
However, $\frac 3 4 \sqrt 2 \simeq 1.0607 >1$, so that this does not prove $P_{n-2}(a) \le P_{n-2}(a_{max})$ for all $a \in S^{n-1}$. However, if $a_1$ satisfies $a_1 \le \sqrt{\sqrt 2 -1} \simeq 0.643$, \eqref{eq16} yields
$$\frac 1 2 P_{n-2}(a) \le \left(n-\frac 3 2+\frac{\sqrt 2 -1} 2\right) \sqrt 2 = (n-2) \sqrt 2 +1 = \frac 1 2 P_{n-2}(a_{max}) . $$
For $a_1 > \frac 1 {\sqrt 2}$, the requirement that $(n-1) \sqrt{1-\frac{1-a_1^2}{n-1}} \frac 1 {a_1} \le (n-2) \sqrt 2 +1 $ is strongest for $n=3$, in which case it means $a_1 \ge \frac 1 {\sqrt{ \sqrt 2 + \frac 1 2 }} \simeq 0.723$. Therefore for any $a \in S^{n-1}$ with
$a_1 \notin (\sqrt{\sqrt 2 -1},\frac 1 {\sqrt{ \sqrt 2 + \frac 1 2 }})$, we have shown $P_{n-2}(a) \le P_{n-2}(a_{max})$. \\
Hence Theorem \ref{th1} is proved also for real scalars, except in the situation that
\begin{equation}\label{eq17}
a_1 \in \left(\sqrt{\sqrt 2 -1},\frac 1 {\sqrt{ \sqrt 2 + \frac 1 2 }}\right) ,
\end{equation}
where the estimate is off by at most $2(\frac 3 4 \sqrt 2 - 1) \simeq 0.121$. This discrepancy occurs since in \eqref{eq15} the extremals for the sum of weights and for the section function $A$ occur for different sequences $a$. This could possibly be avoided, if one could show how the monotonicity of the sequence $(a_k)_{k=1}^n$ affects the size of the integrals $D_k(a)$, but we have not been to find a result of this type. Instead, we will address the case of \eqref{eq17} by a different method in the next section.

\section{Interpolating Ball's function}

To prove Theorem \ref{th1} also for hyperplane sections perpendicular to $a$ with $a_1 \in (\sqrt{\sqrt 2 -1},\frac 1 {\sqrt{ \sqrt 2 + \frac 1 2 }})$, $\KK=\RR$, we will improve the general estimate for $A_{n-1}(a)$ in \eqref{eq16}, by using the improved estimates for Ball's integral function $f$ stated in Propositions \ref{prop6} and \ref{prop7}. The convexity of $f$ allows estimates by interpolation for certain values of $a_1$ and $a_2$ near $\frac 1 {\sqrt 2}$. The technical proof of Proposition \ref{prop6} is given in the Appendix. The proof of Proposition \ref{prop7} is a slight modification of Nazarov, Podkorytov's \cite{NP} proof of Ball's inequality $f(p) \le f(2)=1$ for $p \ge 2$. Recall that Proposition \ref{prop7} states that $f(p) \le \sqrt{\frac 3 \pi} < 1$ for all $p \ge \frac 9 4$. \\

{\bf Proof of Proposition \ref{prop7}.} \\
Let $f:(1,\infty) \to \RR_+$ denote Ball's function, $f(p) := \sqrt{\frac p 2} \frac 2 \pi \ \int_0^\infty |\frac{\sin(t)} t|^p \ dt$. Define
$g, h: [0,\infty) \to \RR_+$ by $g(x):= |\frac{\sin(x)} x|$ and $h(x) := \exp(-\frac{x^2} 6)$, with $g(0)=1$, and let $G, H : (0,1] \to \RR_+$ denote the distribution functions of $g$ and $h$, respectively. We claim that there is $y_0 \in (0,1)$ such that
\begin{equation}\label{eq18}
H(y) \le G(y) \text{ for all } 0 < y < y_0 \quad \text{ and } \quad H(y) \ge G(y) \text{ for all }  y_0 < y < 1 .
\end{equation}
The distribution function lemma in \cite{NP} then implies that the function $\phi :(1,\infty) \to \RR_+$,
$$\phi(p) := \frac 1 {py_0^p} \ \int_0^\infty (h(x)^p - g(x)^p) \ dx $$
is increasing in $p$. Since by Proposition \ref{prop6} for $p_0 := \frac 9 4$
$$\int_0^\infty g(x)^{p_0} = \int_0^\infty \left|\frac{\sin(x)} x\right|^{p_0} dx < \sqrt{\frac 3 \pi} \frac \pi 2 \sqrt{\frac 2 {p_0}} = \sqrt{\frac {2 \pi} 3} =
\int_0^\infty \exp(-\frac 3 8 x^2) dx = \int_0^\infty g(x)^{p_0} dx , $$
we conclude that for all $p \ge \frac 9 4$
$$\int_0^\infty \left|\frac{\sin(x)} x\right|^p \ dx = \int_0^\infty g(x)^p \ dx < \int_0^\infty h(x)^p \ dx = \sqrt{\frac{3 \pi}{2p}} , $$
which is equivalent to $f(p) < \sqrt{\frac 3 \pi}$. For $p=2$,
$$\int_0^\infty g(x)^2 \ dx = \frac \pi 2 > \sqrt{\frac 3 \pi} \frac \pi 2 = \int_0^\infty h(x)^2 \ dx . $$
Therefore there is $q \in (2, \frac 9 4)$ such that
$$0 = \int_0^\infty (h(x)^q - g(x)^q) \ dx = q \int_0^1 y^{q-1} (H(y)-G(y)) \ dy . $$
Hence $H-G$ has at least one zero $y_0 \in (0,1)$. To prove \eqref{eq18}, we will show that $H-G$ has {\em only} one zero. For $m \in \NN$, let
$y_m := \max \{ \ g(x) \ | \ x \in [m \pi,(m+1) \pi] \}$. Since
$$\frac{\sin(x)} x = \prod_{n \in \NN} ( 1 - \frac{x^2}{(n \pi)^2} ) , $$
we have for all $0 < x < \pi$
$$\ln(\frac{\sin(x)} x) = \sum_{n \in \NN} \ln(1 - \frac{x^2}{(n \pi)^2} ) \le - \sum_{n \in \NN} \frac{x^2}{(n \pi)^2} = - \frac{x^2} 6 , $$
i.e. $g(x) = \frac{\sin(x)} x \le \exp(- \frac{x^2} 6) = h(x)$. Therefore $H-G$ is positive in $(y_1,1)$. To show that $H-G$ has only one zero, it
suffices to prove that $(H-G)' > 0$ in $(0,y_1) = \cup_{m \in \NN} [y_{m+1},y_m)$. Since $H'<0, G'<0$, this means that $|\frac {G'}{H'}| >1$ has to be shown there. We have, similarly as in \cite{NP}, $H(y) = h^{-1}(y) = \sqrt{6 \ln (\frac 1 y)}$, $H'(y) = \sqrt{\frac 3 2} \frac 1 {y \sqrt{\ln(\frac 1 y)}}$ and
$$|G'(y)| = \sum_{x>0,\ g(x)=y} \frac 1 {|g'(x)|} . $$
For $y \in (y_{m+1},y_m)$, $g(x)=y$ has one root $x_0$ in $(0,\pi)$ and two roots $x_{j,1}, x_{j,2}$ in $(j \pi, (j+1) \pi)$ for $j = 1, \cdots , m$. Easy estimates show $|g'(x_0)| < \frac 1 2$, $|g'(x_{j,1})|, |g'(x_{j,2})| \le \frac 1 {\pi j}$, $j \in \NN$ so that for all $y \in [y_{m+1},y_m)$ with
$y > y_{m+1} > \frac 1 {\pi (m + \frac 3 2)}$
$$|G'(y)| > 2 (1 + \sum_{j=1}^m  \pi j ) = 2 + \pi m (m+1) , $$
\begin{align*}
\left|\frac{G'(y)}{H'(y)}\right| & > \sqrt{\frac 2 3} (2 + \pi m (m+1)) y \sqrt{\ln(\frac 1 y)} \\
& \ge \sqrt{\frac 2 3} \ \frac{2 + \pi m (m+1)}{\pi (m + \frac 3 2)} \sqrt{\ln(\pi (m +\frac 3 2))} \ge \sqrt{\frac 2 3 \ln(\frac{5 \pi} 2)} > 1 .
\end{align*}
This means that \eqref{eq18} holds and Proposition \ref{prop7} is proven. \hfill $\Box$ \\

Theorem \ref{th1} has been shown for $\KK = \RR$, $n=3$ and for $n \ge 4$ if $a_1 \notin (\sqrt{\sqrt 2 -1},\frac 1 {\sqrt{ \sqrt 2 + \frac 1 2 }})$. We now consider the remaining cases and assume first that $a_2 \ge \frac 2 3$.

\begin{lemma}\label{lem8}
Assume that $a \in S^{n-1}$, $a_1 \in (\sqrt{\sqrt 2 -1},\frac 1 {\sqrt{ \sqrt 2 + \frac 1 2 }})$ and $a_2 \ge \frac 2 3$. Then
$$P_{n-2}(a) \le P_{n-2}(a_{max}) = 2 ( (n-2) \sqrt 2 +1 ) . $$
\end{lemma}

{\bf Proof.} Since $a_n^2+ \cdots + a_3^2 = 1-a_1^2-a_2^2 \le \frac 1 9$, we know that $a_n \le \cdots \le a_3 \le \frac 1 3 < \frac 2 3 \le a_2 \le a_1 \le \frac 1 {\sqrt{ \sqrt 2 + \frac 1 2 }}$. By concavity of $\sqrt{1-x}$, we find similarly as in \eqref{eq16}
\begin{align}\label{eq19}
\frac 1 2 P_{n-2}(a) & \le \sum_{k=2}^n \sqrt{1-a_k^2} \ A_{n-1}(a) \nonumber \\
& \le \left[ \ (n-2) \sqrt{1-\frac{\sum_{k=2}^n a_k^2}{n-2}} + \sqrt{1-a_2^2} \ \right] \ A_{n-1}(a)  \nonumber \\
& = \left[ \ (n-2) \sqrt{1-\frac{1-a_1^2-a_2^2}{n-2}} + \sqrt{1-a_2^2} \ \right] \ A_{n-1}(a) \nonumber \\
& \le \left[ \ (n-2) - \frac{1-a_1^2-a_2^2}{2} + \sqrt{1-a_2^2} \ \right] \ A_{n-1}(a) .
\end{align}
i) Suppose first that $a_1 \le \frac 1 {\sqrt 2}$. Then $ \frac 2 3 \le a_2 \le a_1 \le \frac 1 {\sqrt 2}$, $2 \le a_1^{-2} \le a_2^{-2} \le \frac 9 4$, and by H\"older's inquality with $p_k:=a_k^{-2} \ge \frac 9 4$ for $k \ge 3$ and Proposition \ref{prop7}
\begin{align*}
A_{n-1}(a) &= \frac 2 \pi \int_0^\infty \prod_{k=1}^n \frac{\sin(a_ks)}{a_ks} \ ds
\le \prod_{k=1}^n (\frac 2 \pi \ a_k^{-1} \int_0^\infty \left|\frac{\sin(t)} t\right|^{a_k^{-2}} \ dt)^{a_k^2} \\
& = \left(\prod_{k=1}^n f(a_k^{-2})^{a_k^2}\right) \sqrt2 \le \left(\sum_{k=1}^n a_k^2 f(a_k^{-2})\right) \sqrt 2 \\
& \le [ \ (1-a_1^2-a_2^2) \sqrt{\frac 3 \pi} + a_2^2 f(a_2^{-2}) + a_1^2 f(a_1^{-2}) \ ] \sqrt 2 ,
\end{align*}
where the second inequality follows from the general arithmetic-geometric mean inequality. For $k=1, 2$ write $a_k^{-2}=\lambda_k 2 + (1-\lambda_k) \frac 9 4$, $\lambda_k = 9-4 a_k^{-2}$. Using the convexity of $f$, cf. Proposition \ref{prop6} (c), we find
$$a_k^2 f(a_k^{-2}) \le a_k^2 (\lambda_k f(2) + (1-\lambda_k) f(\frac 9 4)) \le (9 a_k^2 -4) + (4 - 8 a_k^2) \sqrt{\frac 3 \pi} =: \phi_2(a_k) , $$
so that
\begin{align*}
A_{n-1}(a) & \le \left( (1-a_1^2-a_2^2) \sqrt{\frac 3 \pi} + \phi_2(a_2) + \phi_2(a_1) \right) \sqrt 2 \\
& = \left( 9(a_1^2+a_2^2)-8 + 9 (1-a_1^2-a_2^2) \sqrt{\frac 3 \pi}\right) \sqrt 2 =: \psi_1(a_1,a_2) \sqrt 2
\end{align*}
and with \eqref{eq19}
$$\frac 1 2 P_{n-2}(a) \le \left(n-2 - \frac{1-a_1^2-a_2^2}{2} + \sqrt{1-a_2^2} \ \right) \ \psi_1(a_1,a_2) \sqrt 2 =: \gamma(a_1,a_2) \sqrt 2 . $$
Calculation shows that $\frac{\partial \gamma}{\partial a_1} \ge 0$, $\frac{\partial \gamma}{\partial a_2} \ge 0$ for all $a_1,a_2$ in the range considered. Therefore $\gamma$ is increasing in $a_1$ and $a_2$ for all $\frac 2 3 \le a_2 \le a_1 \le \frac 1 {\sqrt 2}$ and $\frac 1 2 P_{n-2}(a)$ is bounded by
$\gamma(\frac 1 {\sqrt 2},\frac 1 {\sqrt 2}) \sqrt 2 = (n-2) \sqrt 2 +1 = \frac 1 2 P_{n-2}(a_{max})$ . \\

ii) Suppose next that $\frac 1 {\sqrt 2} < a_1 \le \frac 1 {\sqrt{\sqrt 2 + \frac 1 2}}$. Write $a_1^{-2} = \lambda \ 2 + (1-\lambda) (\sqrt 2 + \frac 1 2)$. Then by Proposition \ref{prop6}
\begin{align*}
a_1^2 f(a_1^{-2}) & \le a_1^2 ( \lambda f(2) + (1- \lambda) f(\sqrt 2 + \frac 1 2) ) \\
& \le \frac{1-(\sqrt 2 + \frac 1 2) a_1^2 +(2 a_1^2-1) \frac {51}{50}}{\frac 3 2 - \sqrt 2} =: \phi_1(a_1) .
\end{align*}
Using this, we find similarly as in part i)
\begin{align*}
A_{n-1}(a) & \le \left(\sum_{k=1}^n a_k^2 f(a_k^{-2}) \right) \sqrt 2 \\
& \le \left( (1-a_1^2-a_2^2) \sqrt{\frac 3 \pi} + \phi_2(a_2) + \phi_1(a_1) \right) \sqrt 2 =: \psi_2(a_1,a_2) ,
\end{align*}
and with \eqref{eq19}
\begin{align}\label{eq20}
\frac 1 2 P_{n-2}(a) & \le \left((n-2)- \frac{1-a_1^2-a_2^2} 2 + \sqrt{1-a_2^2}  \right) \ \min( \psi_2(a_1,a_2) \sqrt 2 , \frac 1 {a_1} ) \nonumber \\
& =: \min(\gamma_1(a_1,a_2), \gamma_2(a_1,a_2)) ,
\end{align}
where we also used that $A_{n-1}(a) \le \frac 1 {a_1}$ since $a_1 > \frac 1 {\sqrt 2}$. It is easy to see that $\gamma_2$ is decreasing in $a_1$ and in $a_2$ since
$$\frac{\partial \gamma_2}{\partial a_1} = - [ \ (n- \frac 5 2) + \frac 1 2 (a_1^2+a_2^2) + \sqrt{1-a_2^2} \ ] \frac 1 {a_1^2} + 1 < 0 , $$
$$\frac{\partial \gamma_2}{\partial a_2} = \frac{a_2}{a_1} ( 1 - \frac 1 {\sqrt{1-a_2^2}} ) < 0 . $$
A slightly longer calculation and easy estimates show, conversely, that $\gamma_1$ is increasing in $a_1$ and $a_2$. Consider the line
$\frac 1 {\sqrt 2} - a_2 = 8 (a_1 - \frac 1 {\sqrt 2})$ for $\frac 2 3 \le a_2 \le \frac 1 {\sqrt 2} \le a_1 \le \frac 1 {\sqrt{\sqrt 2 + \frac 1 2}}$ (which originates as an approximation of the curve defined by $\gamma_1(a_1,a_2) = \gamma_2(a_1,a_2)$). By \eqref{eq20}
$$\frac 1 2 P_{n-2}(a) \le \max \{ \gamma_1(a_1,a_2), \gamma_2(a_1,a_2) \ | \ \frac 1 {\sqrt 2} - a_2 = 8 (a_1 - \frac 1 {\sqrt 2}) , \frac 2 3 \le a_2 \le \frac 1 {\sqrt 2} \le a_1 \le \frac 1 {\sqrt{\sqrt 2 + \frac 1 2}} \} . $$
One checks that for all $n \ge 4$
$$\frac{\partial \gamma_1}{\partial t}(t,\frac 1 {\sqrt 2}-8 t) = \frac{\partial \gamma_1}{\partial a_1} - 8 \frac{\partial \gamma_1}{\partial a_2} <0 , $$
$$\frac{\partial \gamma_2}{\partial t}(t,\frac 1 {\sqrt 2}-8 t) = \frac{\partial \gamma_2}{\partial a_1} - 8 \frac{\partial \gamma_2}{\partial a_2} <0 . $$
Therefore
$$\gamma_1(a_1,\frac 9 {\sqrt 2} - 8 a_1) \le \gamma_1(\frac 1 {\sqrt 2},\frac 1 {\sqrt 2}) = (n-2) \sqrt 2 + 1 $$
and
$$\gamma_2(a_1,\frac 9 {\sqrt 2} - 8 a_1) \le \gamma_1(\frac 1 {\sqrt 2},\frac 1 {\sqrt 2}) = (n-2) \sqrt 2 + 1 . $$
Hence $P_{n-2}(a) \le 2 ( (n-2) \sqrt 2 +1 ) = P_{n-2}(a_{max})$ for all  $\frac 2 3 \le a_2 \le \frac 1 {\sqrt 2} \le a_1 \le \frac 1 {\sqrt{\sqrt 2 + \frac 1 2}}$. \hfill $\Box$ \\

Next, we consider a similar interpolation scheme if $a_2 < \frac 2 3$.

\begin{lemma}\label{lem9}
Assume that $a \in S^{n-1}$, $a_1 \in (\sqrt{\sqrt 2 -1},\frac 1 {\sqrt 2})$ and $a_2 < \frac 2 3$. Then for all $n \ge 3$
$$P_{n-2}(a) \le P_{n-2}(a_{max}) = 2 ( (n-2) \sqrt 2 +1 ) . $$
If $a_1 \in (\frac 1 {\sqrt 2},\frac 1 {\sqrt{\sqrt 2 + \frac 1 2}})$ and $a_2 < \frac 2 3$, the same holds for all $n \ge 3$, except for possibly $n=5$ or $n=6$.
\end{lemma}

{\bf Proof.} i) Suppose first that $a_1 \le \frac 1 {\sqrt 2}$. Since $a_2 < \frac 2 3$, $f(a_2^{-2}) \le \sqrt{\frac 3 \pi}$ by Proposition \ref{prop7}. If
$\sqrt{\sqrt 2 -1} \le a_1 \le \frac 2 3$, also $f(a_1^{-2}) \le \sqrt{\frac 3 \pi}$. If $\frac 2 3 < a_1 \le \frac 1 {\sqrt 2}$, we again use the convexity of $f$ to get the slightly weaker estimate $a_1^2 f(a_1^{-2}) \le \phi_2(a_1)$, where $\phi_2$ is as in i) of the proof of Lemma \ref{lem8}. This yields
\begin{align*}
A_{n-1}(a) & \le \left( \sum_{k=1}^n a_k^2 f(a_k^{-2}) \right) \sqrt 2 \le \left( a_1^2 f(a_1^{-2}) + (1-a_1^2) \sqrt{\frac 3 \pi} \right) \sqrt 2 \\
& \le \left( (9 a_1^2-4) +(5-9 a_1^2) \sqrt{\frac 3 \pi} \right) \sqrt 2 =: \psi_1(a_1) \sqrt 2 ,
\end{align*}
$$\frac 1 2 P_{n-2}(a) \le \left( (n-1) \sqrt{1-\frac{1-a_1^2}{n-1}} \right) \ \psi_1(a_1) \sqrt 2 =:\gamma(a_1) . $$
As easily seen, $\gamma$ is increasing for $a_1 \le \frac 1 {\sqrt2}$, so that
\begin{align}\label{eq21}
\frac 1 2 P_{n-2}(a) & \le \gamma(\frac 1 {\sqrt 2}) =  (n-1) \sqrt{1 - \frac 1 {2(n-1)}} \ \frac 1 {\sqrt2} \ (1 + \sqrt{\frac 3 \pi}) \nonumber \\
& < (n-2) \sqrt 2 +1 = \frac 1 2 P_{n-2}(a_{max})
\end{align}
for all $n \ge 5$. \\

ii) Assume now that $\frac 1 {\sqrt 2}< a_1 \le \frac 1 {\sqrt{\sqrt 2 + \frac 1 2}}$ and $a_2 < \frac 2 3$. Then again $f(a_2^{-2}) \le \sqrt{\frac 3 \pi}$ and for $f(a_1^{-2})$ we get by interpolation $a_1^2 f(a_1^{-2}) \le \phi_1(a_1)$, where $\phi_1$ is as in part ii) of the proof of Lemma \ref{lem8}. Therefore
$$A_{n-1}(a) \le \left( \sum_{k=1}^n a_k^2 f(a_k^{-2}) \right) \sqrt 2 \le \left( \phi_1(a_1) + (1-a_1^2) \sqrt{\frac 3 \pi} \right) \sqrt 2 =: \psi_2(a_1) \sqrt 2 , $$
\begin{equation}\label{eq22}
\frac 1 2 P_{n-2}(a) \le ( (n-1) \sqrt{ 1 - \frac{1-a_1^2}{n-1}} ) \ \min(\psi_2(a_1) \sqrt 2 , \frac 1 {a_1} ) =: \min ( \gamma_1(a_1), \gamma_2(a_1) ) ,
\end{equation}
where we also used that $A_{n-1}(a) \le \frac 1 {a_1}$ holds. Differentiating $\gamma_1$ and $\gamma_2$, one finds that $\gamma_1'>0>\gamma_2'$ in the range of $a_1$ considered. Therefore, $\gamma_1$ is increasing and $\gamma_2$ is decreasing. We have, independently of $n \in \NN, n \ge 3$, that
$\gamma_1(\bar{a_1}) = \gamma_2(\bar{a_1})$ for $\bar{a_1} \simeq 0.71254$ and $\gamma_1(\bar{a_1}) = \gamma_2(\bar{a_1}) \le (n-2) \sqrt 2 +1$ for all $n \ge 7$. Then $P_{n-2}(a) \le P_{n-2}(a_{max})$. \\

For $n=6$, this estimate is violated by $< 0.006$, for $n=5$ by $< 0.015$. It is correct for $n=5$, $a_1 \notin (0.7095,0.7149)$ and for
$n=6$, $a_1 \notin (0.7115,0.7133)$. \\

For $n=3$ we already proved Theorem \ref{th1}. For $n=4$, the above estimates in i) and ii) yield $P_2(a) \le P_2(a_{max})$ if $a_1 \notin (0.7069,0.7177)$. However, the explicit formulas given in equation \eqref{eq10} for $P_2(a)$ yield $P_2(a) \le P_2(a_{max})$ also for these $a_1$. Only the first or the second case in \eqref{eq10} can occur, since $a_1 < a_2+a_3+a_4$ in our situation. The maximum of the second formula occurs for $a_3=a_4$, with
$P_2(a) < 3.6 < 2( \sqrt 2 +1 )$. The first expression in \eqref{eq10} yields an even smaller maximal value. We do not give the details. In principle, the explicit formulas in \eqref{eq10} could be used to prove $P_2(a) \le P_2(a_{max})$ for all $a \in S^3$, as in the case $n=3$, though this would be more complicated. \\

Replacing \eqref{eq22} by the slightly stronger estimate
$$\frac 1 2 P_{n-2}(a) \le \left( (n-2) \sqrt{1 - \frac{1-a_1^2-a_2^2}{n-2}} + \sqrt{1-a_2^2} \right) \ \min(\psi_2(a_1) \sqrt 2 , \frac 1 {a_1} ) , $$
we get $P_{n-2}(a) \le P_{n-2}(a_{max})$ for all $a_1 \in (\frac 1 {\sqrt 2},\frac 1 {\sqrt{\sqrt 2 + \frac 1 2}})$ and $a_2 \ge 0.5803$ if $n=5$ and $a_2 \ge 0.4952$ if $n=6$. Therefore the only cases left open to prove Theorem \ref{th1} are
\begin{eqnarray}\label{eq23}
\begin{cases}
         & n=5 \quad , \quad a_1 \in (0.7095,0.7149) \quad , \quad a_2 \le 0.5803  \quad\\
         & n=6 \quad , \quad a_1 \in (0.7115,0.7133) \quad , \quad a_2 \le 0.4952  \quad
             \end{cases}     \Bigg\} \ ,
\end{eqnarray}
i.e. when $a_1 \simeq \frac 1 {\sqrt 2}$ and $a_2 < a_1 - \frac 1 8$. This case will be treated by using the following Lemma. \hfill $\Box$ \\

\begin{lemma}\label{lem10}
For $a_1 \in (0.7095,0.7149)$, $\frac 1 {\sqrt{10}} \le a_2 \le 0.5803$ we have
\begin{equation}\label{eq24}
\frac 2 \pi \int_0^\infty  \left|\frac{\sin(a_1s)}{a_1s} \frac{\sin(a_2s)}{a_2s} \right|^{\frac 1 {a_1^2+a_2^2}} \ ds \le 0.985 \sqrt 2 .
\end{equation}
\end{lemma}

We will prove Lemma \ref{lem10} in the Appendix. Using \eqref{eq24}, we finish the proof of Theorem \ref{th1} in the remaining cases \eqref{eq23}: \\

By H\"older's inequality and \eqref{eq24}
\begin{align*}
A_{n-1}(a) & \le (\frac 2 \pi \int_0^\infty \left|\frac{\sin(a_1s)}{a_1s} \frac{\sin(a_2s)}{a_2s} \right|^{\frac 1 {a_1^2+a_2^2}} \ ds )^{a_1^2+a_2^2}
\prod_{j=3}^n (\frac 2 \pi \int_0^\infty |\frac{\sin(a_js)}{a_js}|^{a_j^{-2}} \ ds )^{a_j^2} \\
& \le \left( (a_1^2+a_2^2) \ 0.985 + (1-a_1^2-a_2^2) \sqrt{\frac 3 \pi} \right) \sqrt 2 \le 0.985 \sqrt 2 .
\end{align*}
Therefore
$$\frac 1 2 P_{n-2}(a) \le (n-1) \sqrt{1- \frac{1-a_1^2}{n-1}} \ 0.985 \sqrt 2 . $$
For $n=5$ and $n=6$ and $a_1 \le 0.7149$ this is $<(n-2) \sqrt 2 +1$, so that $P_{n-2}(a) \le P_{n-2}(a_{max})$ also in the cases \eqref{eq23}.
This ends the proof of Theorem \ref{th1}. \hfill $\Box$ \\

As for the lower estimate of $P_{n-2}(a)$, the natural conjecture would be $P_{l(n-2)}(a) \ge P_{l(n-2)}(a_{min}) = 2 \pi^{l-1} (n-1)$, $a_{min}=(1,0,\cdots,0)$, with $l=1$ if $\KK=\RR$ and $l=2$ if $\KK=\CC$. We can only prove a slightly weaker estimate.

\begin{prop}\label{prop11}
For any $a \in \KK^n$ with $|a|=1$
$$P_{n-2}(a) \ge 2 (n-2) \quad , \quad \KK=\RR , $$
$$P_{2(n-2)}(a) \ge 2 \pi (n-2) \quad , \quad \KK=\CC . $$
\end{prop}

{\bf Proof of Proposition \ref{prop11}.} \\
i) We may assume $a \in S^{n-1}$, $a_1 \ge \cdots \ge a_n \ge 0$. In the complex case $\KK=\CC$, by Propositions \ref{prop2}, \ref{prop3} and \ref{prop4}
$$\frac 1 {2 \pi} P_{2(n-2)}(a) \ge \min \left\{ \sum_{k=1}^n (1-a_k^2) C_k \ | \ 0 \le C_k \le A_{2(n-1)}(a) , \sum_{k=1}^n C_k = (n-1) A_{2(n-1)}(a) \right\} . $$
Since $(1-a_k^2)_{k=1}^n$ is increasing in $k$, the sum $\sum_{k=1}^n (1-a_k^2) C_k$ is minimized, if the $C_k$ are decreasing, i.e. for
$C_1= \cdots =C_{n-1}=A_{2(n-1)}(a)$ and $C_1=0$ so that by Hensley \cite{H} and Oleszkiewicz, Pe{\l}czy\'nski \cite{OP}
\begin{align*}
\frac 1 {2\pi} P_{2(n-2)}(a) & \ge \sum_{k=1}^{n-1} (1-a_k^2) A_{2(n-1)}(a) = (n-1 -(1-a_1^2)) A_{2(n-1)}(a) \\
& \ge (n-2) A_{2(n-1)}(a) \ge (n-2) .
\end{align*}
ii) Similarly, we find in the real case $\KK=\RR$, using Hensley's lower estimate \cite{H} for the parallel section function $A_{n-1}$
$$\frac 1 2 P_{n-2}(a) \ge \sum_{k=1}^{n-1} \sqrt{1-a_k^2} \ A_{n-1}(a) \ge \sum_{k=1}^{n-1} \sqrt{1-a_k^2} . $$
Now $\phi(x) = \sqrt{1-x}$ is concave and decreasing on $[0,1]$. Therefore for any $x_2<y_2<y_1<x_1$ with $x_2^2+x_1^2 = y_2^2+y_1^2$ we have that
$\phi(y_2^2)+\phi(y_1^2) \ge \phi(x_2^2)+\phi(x_1^2)$, i.e. the sum gets smaller by moving all coordinates towards 0 and 1. Hence
$$\frac 1 2 P_{n-2}(a) \ge (n-2)  + \sqrt{1- \sum_{j=1}^{n-1} a_j^2} \ge (n-2) . $$   \hfill $\Box$ \\

{\bf Remarks.}\\
(a) One possibility to improve the lower estimate in Proposition \ref{prop11} would be to understand how the monotonicity properties of the sequence $a=(a_k)_{k=1}^n$ affect the size of the integrals $D_k(a)$, since e.g. in the real case
$$P_{n-1}(a)  = \sum_{k=1}^n \sqrt{1-a_k^2} \ D_k(a) . $$ \\

(b) Numerical estimates of Ball's integral function $f$, $f(p):= \sqrt{\frac p 2} \frac 2 \pi \int_0^\infty |\frac{\sin(t)} t|^p \ dt$ indicate that $f(p)=\sqrt{\frac 3 \pi}$ for some $p_1$ with $2.165 < p_1 < 2.166$, that $f$ attains its minimum in $p_2$ with $3.36 < p_2 < 3.37$ and that $f$ is convex for $1 < p < p_0$ with $4.46 < p_0 < 4.47$. The behavior of $f$ near $\infty$ is well-understood: by a result of Kerman, Ol'hava and Spektor \cite{KOS}
$$f(p) = \sqrt{\frac 3 \pi} \left(1 - \frac 3 {20} \frac 1 p - \frac {13}{1120} \frac 1 {p^2} + O(\frac 1 {p^3}) \right) . $$

\section{An application  of the Busemann-Petty type}\label{BP}
In this section we apply the result of Theorem \ref{th1} to the surface area version of the Busemann-Petty
problem described in the Introduction.

\begin{theorem} For each $n\ge 14,$ there exist origin-symmetric convex bodies $K,L$ in $\RR^n$
such that for all $a\in S^{n-1}$
$$vol_{n-2}(\partial K\cap a^\bot)\le vol_{n-2}(\partial L\cap a^\bot)$$
but
$$vol_{n-1}(\partial K)> vol_{n-1}(\partial L).$$
\end{theorem}

{\bf Proof.} Let $K=B_\infty^n$ be the unit cube in $\RR^n.$ Let $L$ be the Euclidean ball of radius $r$ in $\RR^n$
so that the perimeters of hyperplane sections of $L$ are all equal to the maximal perimeter of sections
of $K.$ Namely, for any $a\in S^{n-1}$

$$vol_{n-2}(\partial K\cap a^\bot ) \le vol_{n-2}(\partial K\cap a_{max}^\bot) = 2 ((n-2) \sqrt 2 +1) = vol_{n-2}(r S^{n-2}) = r^{n-2} \ \frac{2 \pi^{(n-1)/2}}{\Gamma(\frac {n-1} 2)} \ , $$
i.e. $$r = \frac{ [((n-2) \sqrt 2 +1)\Gamma(\frac {n-1} 2)]^{\frac 1 {n-2}} }{\pi^{(n-1)/(2(n-2))}}.$$

The desired inequality for the surface areas of $K$ and $L$ happens when
$$vol_{n-1}(\partial B_\infty^n) = 2 n > vol_{n-1}(r S^{n-1}) = r^{n-1} \ \frac{2 \pi^{n/2}}{\Gamma(\frac n 2)} \ .$$
The latter is equivalent to
$$1 > \frac{\pi^{n/2}}{n \Gamma(\frac n 2)} r^{n-1} =  \frac{\pi^{n/2}}{n \Gamma(\frac n 2)}   \frac{ [((n-2) \sqrt 2 +1)\Gamma(\frac {n-1} 2)]^{\frac {n-1}{n-2}} }{\pi^{(n-1)^2/(2(n-2))}}$$$$   = \frac 1 {n \Gamma(\frac n 2)}   \frac{ [((n-2) \sqrt 2 +1)\Gamma(\frac {n-1} 2)]^{\frac {n-1}{n-2}} }{\pi^{1/(2(n-2))}} =: BP(n) \ . $$
Then $BP$ is decreasing in $n$, with $BP(x_0)=1$ for $x_0 \simeq 13.70$, so $BP(n)<1$ for all $n \ge 14$. \qed\\

A similar argument can be made in the complex case when there are similar counterexamples for all $n \ge 11$.

\section{Appendix}

In the Appendix, we present the technical proofs of Proposition \ref{prop6} and of Lemma \ref{lem10}. \\

{\bf Proof of Proposition \ref{prop6} (a).} \\
The fact that $\lim_{p \to \infty} f(p) = \sqrt{\frac 3 \pi}$ is well-known \cite{KOS}, following from $\frac{\sin(x)} x \le \exp(-\frac {x^2} 6)$ for $0 \le x \le \pi$ and $\sqrt{\frac p 2} \frac 2 \pi \int_0^\infty \exp(-\frac{x^2 p} 6) \ dx = \sqrt{\frac 3 \pi}$. Now let $p_0:= \frac 9 4$. Since
$$\frac{\sin(x)} x = \prod_{n \in \NN} \left(1- \frac{x^2}{(n\pi)^2}\right) , $$
and $\ln(1-y) \le -y - \frac 1 2 y^2$ for $0 \le y < 1$, we find for $0 < x < \pi$
$$\ln\left(\frac{\sin(x)} x\right) = \sum_{n \in \NN} \ln(1- \frac{x^2}{(n\pi)^2}) \le - \sum_{n \in \NN} (\frac x {n \pi})^2 - \frac 1 2 \sum_{n \in \NN} (\frac x {n \pi})^4 = - \frac {x^2} 6 - \frac{x^4} {180} , $$
$$\frac{\sin(x)} x \le \exp(-\frac{x^2} 6 - \frac{x^4}{180}) \quad , \quad x \in (0,\pi) . $$
This implies
$$\left( \frac{\sin(x)} x \right)^{p_0} \le \left( \frac{\sin(x)} x \right)^2 \ \exp(- \frac{x^2}{24} - \frac{x^4}{720}) \le \left( \frac{\sin(x)} x \right)^2 \ (1- \frac {x^2}{24}) \;
 , \; x \in (0,\pi) , $$
\begin{align*}
I_0 & := \sqrt{\frac {p_0} 2} \frac 2 \pi \int_0^\pi \left|\frac{\sin(x)} x\right|^{p_0} \ dx \le \sqrt{\frac {p_0} 2} \frac 2 \pi \int_0^\pi \left|\frac{\sin(x)} x\right|^2 \ (1- \frac {x^2}{24}) \ dx \\
& = \sqrt{\frac {p_0} 2} \frac 2 \pi ( Si(2 \pi) - \frac \pi {48} ) \le 0.91340 .
\end{align*}
Here $Si$ denotes the sine integral function.
For $x \in (\pi, 2 \pi)$, $|\frac{\sin(x)} x|^{1/4} \le |\frac{\sin(x_0)} {x_0}|^{1/4} \le 0.683$ where $x_0 \simeq 4.493$. Hence
\begin{align*}
I_1 &:= \sqrt{\frac {p_0} 2} \frac 2 \pi \int_\pi^{2 \pi} \left|\frac{\sin(x)} x\right|^{p_0} \ dx
\le 0.683 \sqrt{\frac {p_0} 2} \frac 2 \pi \int_\pi^{2 \pi} \left|\frac{\sin(x)} x\right|^2 \ dx \\
& = 0.683 \sqrt{\frac {p_0} 2} \frac 2 \pi (Si(4 \pi) -Si(2 \pi)) \le 0.03414.
\end{align*}
For $x \in (k \pi, (k+1) \pi)$, $|\frac{\sin(x)} x|^{p_0} \le (\frac 1 {k \pi})^{1/4} \ (\frac{\sin(x)} x)^2$ and
$$I_k := \sqrt{\frac {p_0} 2} \frac 2 \pi \int_{k \pi}^{(k+1) \pi} \left|\frac{\sin(x)} x\right|^{p_0} \ dx
\le \sqrt{\frac {p_0} 2} \frac 2 \pi (\frac 1 {k \pi})^{1/4} \ (Si((2k+1) \pi) - Si(2k \pi)) . $$
Also,
\begin{align*}
J_k & := \sqrt{\frac {p_0} 2} \frac 2 \pi \int_{k \pi}^\infty \left|\frac{\sin(x)} x\right|^{p_0} \ dx
\le \sqrt{\frac {p_0} 2} \frac 2 \pi  (\frac 1 {k \pi})^{1/4} \int_{k \pi}^\infty \left|\frac{\sin(x)} x\right|^2 \ dx  \\
& = \sqrt{\frac {p_0} 2} \frac 2 \pi  (\frac 1 {k \pi})^{1/4} ( \frac \pi 2 - Si(2k \pi) ) .
\end{align*}
Calculation then shows that $f(p_0) = \sum_{k=0}^5 I_k + J_6 \le 0.977 < \sqrt{\frac 3 \pi} $. \hfill $\Box$ \\

{\bf Proof of Proposition \ref{prop6} (b).} \\
Now let $p_0:=\sqrt 2 + \frac 1 2 \simeq 1.9142 < 2$. The claim is that $f(p_0) < \frac {51}{50}$. We note that for $p<2$, $f(p)>1$. For $0<x<\frac \pi 2$ and $p>1$, we have, similarly as above,
$$\left(\frac{\sin(x)} x \right)^p \le \exp(-\frac{x^2} 6) \ (1-\frac p {180} x^4) , $$
yielding
\begin{align*}
I_{0,1} & :=  \sqrt{\frac {p_0} 2} \frac 2 \pi \int_0^{\frac \pi 2} \left|\frac{\sin(x)} x\right|^{p_0} \ dx
\le  \sqrt{\frac {p_0} 2} \frac 2 \pi \int_0^{\frac \pi 2} \exp(-\frac{x^2} 6) \ (1-\frac {p_0} {180} x^4) \ dx \\
& = \sqrt{\frac {p_0} 2} \frac 2 \pi \left[\frac \pi {480 p_0} (36+\pi^2 p_0) \exp(-\frac{\pi^2 p_0}{24}) + \frac{\sqrt{6 \pi}}{p_0^{3/2}}(\frac{p_0} 2 - \frac 3 {40}) \ \text{erf}(\pi \sqrt{\frac{p_0}{24}})\right] \le 0.76509 ,
\end{align*}
where erf denotes the standard error function, erf$(x):= \frac 2 {\sqrt \pi} \int_0^\infty \exp(-t^2) \ dt$. For $x \in (\frac \pi 2 , \pi)$, Taylor expansion at $x=\pi$ yields an approximation
$$\left(\frac{\sin(x)} x\right)^2 \le (1- \frac x \pi)^2 \sum_{j=0}^6 c_j (1- \frac x \pi)^j , $$
where $c_0=1$, $c_1=2$, $c_2=-(\frac{\pi^2} 3 -3)$ etc. Taking the $\frac{p_0} 2$-th power of this gives an estimate for
$\left(\frac{\sin(x)} x\right)^{p_0}$, $x \in (\frac{\pi} 2 , \pi)$ of the form
$$(\frac{\sin(x)} x)^{p_0} \le (1- \frac x \pi)^{p_0} \sum_{j=0}^6 d_j (1- \frac x \pi)^j \ , $$
$d_0=1$, $d_1=1.9142$ etc., where the right hand side may be integrated exactly, giving
$$I_{0,2} := \sqrt{\frac {p_0} 2} \frac 2 \pi \int_{\frac \pi 2}^{\pi} \left|\frac{\sin(x)} x\right|^{p_0} \ dx \le 0.13531 . $$
In $(\pi, \infty)$, we use approximations of $(\sin(x))^2$ and $|\sin(x)|^{p_0}$ instead of those for $|\frac{\sin(x)} x|^{p_0}$. For $x \in (0,2 \pi)$, we have the alternating series estimate
$$(\sin(x))^2 \le \sum_{j=0}^4 (-1)^j c_j (x-\frac \pi 2)^{2j} \ , \ c_0=c_1=1, c_2= \frac 1 3, c_3=\frac 2 {45}, c_4 = \frac 1 {315} . $$
Taking the $\frac{p_0} 2$-th power of this gives an estimate for
$$(\sin(x))^{p_0} \le \sum_{j=0}^4 (-1)^j d_j (x-\frac \pi 2)^{2j} =: l(x) \ , \ d_0=1, d_1= \frac{p_0} 2, d_2= \frac{p_0^2} 8 -\frac{p_0}{12} \text {  etc.} $$
Using this, we find for $l \in \NN$
\begin{align*}
\sqrt{\frac {p_0} 2} \frac 2 \pi \int_\pi^{(l+1) \pi} \left|\frac{\sin(x)} x\right|^{p_0} \ dx & =
\sqrt{\frac {p_0} 2} \frac 2 \pi \sum_{k=1}^l \int_0^{\pi} \frac{\sin(x)^{p_0}} {(x+k \pi)^{p_0}} \ dx \\
& =: \sum_{k=1}^l I_k  \le \sqrt{\frac {p_0} 2} \frac 2 \pi \sum_{k=1}^l \int_0^{\pi} \frac{l(x)} {(x+k \pi)^{p_0}} \ dx .
\end{align*}
The integrals $\int_0^{\pi} \frac{l(x)} {(x+k \pi)^{p_0}} \ dx$ can be integrated exactly. One finds, choosing $l=5$,
$$(I_{0,1} + I_{0,2}) + (\sum_{k=1}^5 I_k) \le 0.90040 + 0.09383 = 0.99423 . $$
The remaining integral over $(6 \pi, \infty)$ is estimated by
\begin{align*}
\sqrt{\frac {p_0} 2} \frac 2 \pi \int_{6 \pi}^\infty  \left|\frac{\sin(x)} x\right|^{p_0} \ dx
& \le  \sqrt{\frac {p_0} 2} \frac 2 \pi  \sum_{k=6}^\infty \frac 1 {(k \pi)^{p_0}} \int_0^\pi l(x) dx \\
& \le  \sqrt{\frac {p_0} 2} \frac 2 \pi \frac 1 {\pi^{p_0}} \ ( \zeta(p_0) - \sum_{k=1}^5 \frac 1 {k^{p_0}} ) \ 1.6 \le 0.02567 .
\end{align*}
Hence
$f(p_0) \le 0.99423 + 0.02567 =1.0199 < \frac {51}{50}$ . \hfill $\Box$ \\

{\bf Proof of Proposition \ref{prop6} (c).} \\
(a) We claim that $f|_{[\sqrt2 + \frac 1 2, \frac 9 4]}$ is convex. Put $g(p,x):= \sqrt p \ |\frac{\sin(x)} x|^p$ and \\
$h(p,x):= 4 p^2 (\ln |\frac{\sin(x)} x|)^2 + 4 p (\ln |\frac{\sin(x)} x|) - 1$. Then
$$\frac{\partial^2 g}{\partial p^2}(p,x) = h(p,x) \ \frac 1 {4 p^{\frac 3 2}} \ \left|\frac{\sin(x)} x\right|^p , x>0 . $$
We will show that for all $p \in [\sqrt 2 + \frac 1 2, \frac 9 4]$
\begin{equation}\label{eq25}
\int_0^\infty h(p,x) \ \left|\frac{\sin(x)} x\right|^p \ dx \ge \frac 1 5 >0 .
\end{equation}
This yields that $f''(p) >0$, so that $f$ is convex in this interval. We remark that $f$ is negative in $0 < x < x_p$ and non-negative for $x>x_p$, where $x_p \simeq 1.8$ for $p$ in the given range: More precisely, note that
$$4 p^2 (\ln y)^2 + 4 p (\ln y) -1 = 0 \; \Leftrightarrow \; \ln y = - \frac{1 \pm \sqrt 2}{2 p} \; \Leftrightarrow \; y = \exp(- \frac{1 \pm \sqrt 2}{2 p})  \ . $$
For $y = |\frac{\sin(x)} x| < 1$, we need the plus sign and $y_p = \exp(-\frac{1 + \sqrt 2}{2 p}) = \frac{\sin(x_p)}{x_p}$ yields a unique solution $x_p$, decreasing with $p$,
\begin{eqnarray*}
x_p = \begin{cases} & 1.8205 \quad p = \sqrt 2 + 1/2 \\
                    & 1.7863 \quad p = 2 \\
                    & 1.6965 \quad p = 9/4
                    \end{cases}  \ ,
\end{eqnarray*}
\begin{eqnarray*}
y_p = \begin{cases} & 0.5323 \quad p = \sqrt 2 + 1/2 \\
                    & 0.5469 \quad p = 2 \\
                    & 0.5848 \quad p = 9/4
                    \end{cases}  \ .
\end{eqnarray*}
We estimate $\frac{\partial^2 g}{\partial p^2} (p,x)$ from below on the intervals $(0,x_p)$, $(x_p,\pi)$ and $(\pi,2 \pi)$ and integrate that. \\

(b) For all $0 < x < x_p$, $h(p,x) < 0$. For these we have to estimate $(\frac{\sin(x)} x)^p$ from above to estimate the integrand in \eqref{eq25} from below. We again use
$$\left(\frac{\sin(x)} x\right)^p \le \exp(-\frac{x^2 p} 6) \ (1- \frac p {180} x^4) $$
as in the previous proof. Replacing $\ln(\frac{\sin(x)} x)$ in $h(p,x)$ by $-\frac{x^2} 6 - \frac{x^4}{180}$ increases $h(p,x)$ by at most $10^{-3}$, and only for $x \le 1$, so that
\begin{align*}
\int_0^{x_p} h(p,x) & \ \left(\frac{\sin(x)} x\right)^p \ dx \\
 & \ge \int_0^{x_p} \exp(-\frac{x^2 p} 6) \ (1- \frac p {180} x^4) \ (4 p^2 (\frac{x^2} 6 + \frac{x^4}{180})^2 + 4 p  (-\frac{x^2} 6 - \frac{x^4}{180}) -1) \ dx - 10^{-3} \\
 & \ge \int_0^{x_p} \exp(-\frac{x^2 p} 6) \ (-1 - \frac 2 3 p x^2 +(\frac{p^2} 9 - \frac p {60}) x^4 +\frac {p^2}{900} x^6) \ dx - 10^{-3} \\
 & = - \frac 1 3 x_p \ (3 + p x_p^2 + \frac p {100} x_p^4) \ \exp(-\frac p 6 x_p^2) - 10^{-3} =: \gamma_1(p) \ .
\end{align*}
The second inequality follows by expanding the product of both polynomials and easy lower estimates. Actually, the leading term in $x^6$ is
$\frac {p^2}{90} x^6$, but the lesser value $\frac {p^2}{900} x^6$ was chosen to allow for an exact integration without error functions. We note that
$p x_p^2$ and $\gamma_1(p)$ are both increasing in $p$, $\gamma_1$ with negative values,
\begin{eqnarray*}
\gamma_1(p) = \begin{cases} & -2.015 \quad p = \sqrt 2 + 1/2 \\
                            & -1.971 \quad p = 2 \\
                            & -1.858 \quad p = 9/4
                            \end{cases}  \ .
\end{eqnarray*}

(c) For $x \in (x_p,\pi)$, $h(p,x) >0$. For $y = \frac{\sin(x)} x$, $\frac{dy}{dx} = \frac{\cos(x)} x - \frac{\sin(x)}{x^2}$. We have
$-0.4362 \le \frac{dy}{dx} \le \frac 1 \pi \simeq 0.3183$ for all $x \in (1.6965,\pi)$. Therefore
\begin{align*}
\int_{x_p}^\pi h(p,x) \ \left|\frac{\sin(x)} x\right|^p \ dx & \ge \int_{x_p}^\pi h(p,x) \ \left|\frac{\sin(x)} x\right|^p \ \left|\frac{d(\frac{\sin(x)} x)}{dx}\right| \ \frac 1 {0.4362} \ dx \\
& \ge 2.2926 \ \int_0^{y_p} (4 p^2 (\ln y)^2 + 4 p (\ln y) -1) \ y^p \ dy =: I \ ,
\end{align*}
substituting $y=\frac{\sin(x)} x$, which maps $(x_p,\pi)$ bijectively onto $(0,y_p)$. The last integral can be calculated explicitly,
$$I =  2.2926 \ \left( \frac{4 p^2}{1+p} \ (\ln y_p)^2 -\frac{4 p (p-1)}{(1+p)^2} \ (\ln y_p) + \frac{3 p^2-6 p -1}{(1+p)^3} \right) \ y_p^{p+1} =: \gamma_2(p) \ . $$
The function $\gamma_2$ is increasing in $p$, too, with positive values. One has e.g.
\begin{eqnarray*}
\gamma_2(p) = \begin{cases} & 0.898 \quad p = \sqrt 2 + 1/2 \\
                            & 0.916 \quad p = 2 \\
                            & 0.956 \quad p = 9/4
                            \end{cases}  \ .
\end{eqnarray*}

(d) For $x \in(\pi,2 \pi)$, $h(p,x)$ attains large values: if $p \le 2$, $h(p,x) \ge 21$ and if $p>2$, even $h(p,x) > 24$. We find for $p \le 2$
\begin{align*}
\int_\pi^{2 \pi} h(p,x) \ \left|\frac{\sin(x)} x\right|^p \ dx & \ge 21 \int_\pi^{2 \pi} \left|\frac{\sin(x)} x\right|^p \ dx \ge 21 \int_\pi^{2 \pi} \left|\frac{\sin(x)} x\right|^2 \ dx \\
& = 21 ( Si(4 \pi)- Si(2 \pi) ) \simeq 1.554 =: \gamma_3(p) ,
\end{align*}
and for $p > 2$, using H\"older's inequality,
\begin{align*}
\int_\pi^{2 \pi} h(p,x) & \ \left|\frac{\sin(x)} x\right|^p \ dx  \ge 24 \int_\pi^{2 \pi} \left|\frac{\sin(x)} x\right|^p \ dx \\
& \ge 24  \frac 1 {\pi^{p/2-1}} (\int_\pi^{2 \pi} \left|\frac{\sin(x)} x\right|^2 \ dx)^{p/2} = 24 \pi (\frac 1 \pi [Si(4 \pi) - Si(2 \pi)])^{p/2} =: \gamma_3(p) \ .
\end{align*}
The function $\gamma_3$ is decreasing in $[2,\frac 9 4]$, $\gamma_3(2) \simeq 1.776$, $\gamma_3(\frac 9 4) \simeq 1.112$. \\

We conclude from the estimates in (b), (c) and (d) that
$$\int_0^\infty h(p,x) \ \left|\frac{\sin(x)} x\right|^p \ dx \ge \int_0^{2 \pi} h(p,x) \ \left|\frac{\sin(x)} x\right|^p \ dx \ge \gamma_1(p) + \gamma_2(p) + \gamma_3(p) =: \gamma(p) \ . $$
This way, we find that $\gamma(p) \ge 0.43$ for $\sqrt 2 + \frac 1 2 \le p \le 2$ and $\gamma(p) \ge 0.21$ for $2 < p \le \frac 9 4$. Hence
$f''|_{[\sqrt 2 + \frac 1 2, \frac 9 4]} >0$, and $f$ is convex there. \\
Since $f(\sqrt 2+\frac 1 2) > 1 > f(\frac 9 4)$ and $f$ is convex, $f$ is decreasing in $[\sqrt 2 + \frac 1 2, \frac 9 4]$.  \hfill $\Box$ \\

It remains to prove Lemma \ref{lem10} which is only used for $\KK=\RR$, $n=5$ or $n=6$ and $a_1$ in a very small interval near $\frac 1 {\sqrt 2}$. Therefore we only outline the essential parts of the proof. Basically, the integral in \eqref{eq24} is strictly less than $\sqrt 2$ since $a_1$ and $a_2$ deviate by at least $\frac 1 8$, and so there is some cancelation in the product $\sin(a_1s) \sin(a_2s)$. \\

{\bf Proof of Lemma \ref{lem10}.} We have to estimate with $d := \frac{a_2}{a_1}$
$$\frac 2 \pi \int_0^\infty  \left|\frac{\sin(a_1r)}{a_1r} \frac{\sin(a_2r)}{a_2r} \right|^{\frac 1 {a_1^2+a_2^2}} \ dr
= \sqrt 2 \frac 1 {\sqrt 2 a_1} \frac 2 \pi \int_0^\infty  \left|\frac{\sin(s)}{s} \frac{\sin(ds)}{ds} \right|^{\frac 1 {a_1^2(1+d^2)}} \ ds  =: \sqrt 2 J \ . $$
Note that, by assumption, $\sqrt 2 a_1 \simeq 1$, $a_1^{-2} \simeq 2$. We have to show $J \le 0.985$. \\

a) For $0 < s < \pi$, $\ln (\frac{\sin(x)} x) \le -\frac{s^2} 6 - \frac{s^4}{180} - O(s^6)$, with only negative terms in the series. This yields
$$\left|\frac{\sin(s)}{s} \frac{\sin(ds)}{ds} \right|^{\frac 1 {a_1^2(1+d^2)}} \le \exp(-\frac{s^2}{6 a_1^2}) \ (1- \frac{s^4}{180 \ a_1^2} \frac{1+d^4}{1+d^2}) \ ,$$
Integration of the right side over $[0,\pi]$ yields
\begin{align*}
I_0 & := \frac 1 {\sqrt 2 a_1} \frac 2 \pi \int_0^\pi \left|\frac{\sin(s)}{s} \frac{\sin(ds)}{ds} \right|^{\frac 1 {a_1^2(1+d^2)}} \ ds \\
& \le \sqrt{\frac 3 \pi} \text{erf}(\frac {\pi}{\sqrt 6 a_1}) (1-\frac 3{20}a_1^2\frac{1+d^4}{1+d^2}) + \frac{\sqrt 2}{60} \exp(-\frac{\pi^2}{6 a_1^2})
\frac{\pi^2+9 a_1^2}{a_1} \frac{1+d^4}{1+d^2} = \gamma_1(a_1,d) \ .
\end{align*}

b) For $\pi \le s \le 2 \pi$, we use H\"older's inequality with exponent $p = 2 a_1^2 (1+d^2)$,
\begin{align*}
I_1 & := \frac 1 {\sqrt 2 a_1} \frac 2 \pi \int_\pi^{2 \pi} \left|\frac{\sin(s)}{s} \frac{\sin(ds)}{ds} \right|^{\frac 1 {a_1^2(1+d^2)}} \ ds \\
& \le \frac 1 {\sqrt 2 a_1} \left(\frac 2 \pi \int_\pi^{2 \pi} \left|\frac{\sin(s)}{s} \frac{\sin(ds)}{ds} \right|^2 \ ds\right)^{\frac 1 {2 a_1^2 (1+d^2)}} \
(\frac 2 \pi \pi)^{1-\frac 1 {2 a_1^2 (1+d^2)}} \\
& = \frac {\sqrt 2}{a_1} [ \ (Si(2 \pi)-Si(4 \pi)) \frac 1 {3 \pi d^2} + (Si(2 \pi d)-Si(4 \pi d)) \frac d {3 \pi} \\
& \quad \quad - \sum_{\pm} (Si(2 \pi (1 \pm d)) - Si(4 \pi (1 \pm d))) \frac{(1 \pm d)^3}{6 \pi d^2} - \sin(\pi d)^2 \ \frac{\cos(2 \pi d)}{3 \pi^2 d^2} \ ]^{\frac 1 {2 a_1^2 (1+d^2)}} =:\gamma_2(a_1,d) \ .
\end{align*}
The function $\gamma_2$ is increasing in $a_1$, since $\frac{A^{1/a_1^2}}{a_1}$ is increasing in $a_1$ if $A < \exp(-\frac{a_1^2} 2) \simeq \exp(-\frac 1 4)$ which is satisfied if $A$ is the appropriate power of the integral. \\

c) For $2 \pi \le s < \infty$, we again use H\"older's inequality with exponent $p = 2 a_1^2 (1+d^2)$, but differently
\begin{align*}
I_2 & := \frac 1 {\sqrt 2 a_1} \frac 2 \pi \int_{2 \pi}^\infty \left|\frac{\sin(s)}{s} \frac{\sin(ds)}{ds} \right|^{\frac 1 {a_1^2(1+d^2)}} \ ds \\
& \le \frac 1 {\sqrt 2 a_1} \frac 2 \pi \left(\int_{2 \pi}^\infty \left|\frac{\sin(s)}{s}\right|^2 \ ds\right)^{\frac 1 {2 a_1^2(1+d^2)}} \
\left(\int_{2 \pi}^\infty |\frac 1 {ds}|^{\frac 2 {a_1^2 (1+d^2) -1}} \ ds\right)^{1-\frac 1 {2 a_1^2(1+d^2)}} \ .
\end{align*}
Here we estimated $|\sin(d s)|$ by 1. This yields
$$I_2 \le \frac {\sqrt 2} {\pi a_1} \left(\frac \pi 2 -Si(4 \pi)\right)^{\frac 1 {2 a_1^2(1+d^2)}} \
\left[ \ (\frac 1 {2 \pi d})^{\frac 2 {2 a_1^2 (1+d^2) -1}} \ \frac{2 \pi (2 a_1^2 (1+d^2) -1)}{3 - 2 a_1^2 (1+d^2)} \ \right]^{1-\frac 1 {2 a_1^2(1+d^2)}} =: \gamma_3(a_1,d) \ . $$

In conclusion, $J = I_0+I_1+I_2 \le (\gamma_1+\gamma_2+\gamma_3)(a_1,d) =: \gamma(a_1,d)$, which, for the range of $a_1$, $a_2$ and $d = \frac{a_2}{a_1}$ considered, is bounded by $J \le 0.985$. The function $\gamma$ essentially does not depend on $a_1$ since it is considered only in a tiny interval. Nevertheless, its maximum (for the relevant values of $d$) is attained for the maximal choice $a_1 = 0.7149$. As a function of $d$, the maximum is attained for the maximal value of $d$, which is $d = \frac{0.5803}{0.7095} \le 0.818$, i.e. $J \le \gamma(0.7149,0.818) \le 0.985$. This is not surprising, since the cancelation effect in $\sin(s) \sin(d s)$ is minimal if $d$ is maximal. Obviously, the main contribution to $J$ comes from $I_0$, which is slightly larger than 0.92 whereas $I_1$ and $I_2$ each contribute about 0.03.    \hfill $\Box$ \\

\vspace*{1cm}

\noindent Mathematisches Seminar\\
Universit\"at Kiel\\
24098 Kiel, Germany\\
hkoenig@math.uni-kiel.de\\
\\

\noindent Department of Mathematics\\
University of Missouri\\
Columbia, MO 65211, USA\\
koldobskiya@missouri.edu


\vspace{1cm}



\begin{thebibliography}{999}
\bibitem[B1]{B1} K. Ball; Cube slicing in $\RR^n$, Proc. Amer. Math. Soc. 97 (1986), 465-473.
\bibitem[B2]{B2}  K. Ball; {Some remarks on the geometry of convex sets},
Geometric aspects of functional analysis (1986/87), Lecture Notes in Math. 1317,
Springer-Verlag, Berlin-Heidelberg-New York, 1988, 224--231.
\bibitem[GR]{GR} I, S. Gradstein, I. M. Ryshik; Table of Series, Products and Integrals, VEB Deutscher Verlag der Wiss., Berlin, 1957.
\bibitem[H]{H} D. Hensley; Slicing the cube in $\RR^n$ and probability, Proc. Amer. Math. Soc. 73 (1979), 95-100.
\bibitem[KOS]{KOS} R. Kerman, R. Ol'hava, S. Spektor; An asymptotically sharp form of Ball's integral inequality, Proc. Amer. Math. Soc. 143 (2015), 3839-3846.
\bibitem[K]{K} A. Koldobsky; Fourier analysis in convex geomety, Math. Surv. and Monogr. 116, Amer. Math. Soc., 2005.
\bibitem[KK1]{KK1} H. K\"onig, A. Koldobsky; Volumes of low-dimensional slabs and sections in the cube, Adv. in Appl. Math. 47 (2011), 894-907.
\bibitem[KK2]{KK2} H. K\"onig, A. Koldobsky; Minimal volume of slabs in the complex cube, Proc. Amer. Math. Soc. 140 (2012), 1709-1717.
\bibitem[NP]{NP} F. L. Nazarov, A. N. Podkorytov; Ball, Haagerup and distribution functions, Complex analysis, operators, and related topics, Oper. Th. Adv. Appl. 113, Birkhaeuser, 2000, 247-267.
\bibitem[OP]{OP} K. Oleszkiewicz, A. Pe{\l}czy\'nski; Polydisc slicing in $\CC^n$, Studia Math. 142 (2000), 281-294.
\bibitem[P]{P} A. Pe{\l}czy\'nski; Personal communication, 2007.
\bibitem[Po]{Po} G. P\'olya; Berechnung eines bestimmten Integrals, Math. Ann. 74 (1913) 204-212.
\bibitem[W]{W} G. N. Watson, A treatise on the theory of Bessel functions, Cambr. Univ. Press, 1995.
\end{thebibliography}
\end{document}